\documentclass{amsart} 

\usepackage[all]{xy}
\usepackage{amssymb} 
\usepackage{amsmath} 
\usepackage{amscd} 
\usepackage{latexsym} 
\usepackage{amsthm} 
\usepackage{verbatim} 
\usepackage{graphics} 
\usepackage{epsfig}
 
\let\cal\mathcal

\def\Dscr{{\cal D}} 
\def\Escr{{\cal E}} 
\def\Fscr{{\cal F}}

\def\Lscr{{\cal L}} 
 
\def\Nscr{{\cal N}} 
\def\Oscr{{\cal O}} 
 
\def\Qscr{{\cal Q}}

\let\blb\mathbb 
 
\def\CC{{\blb C}}

\def \PP{{\blb P}} 
\def \ZZ{{\blb Z}}

\def\1{{\mathbb{1}}}

\def\codim{\operatorname {codim}} 
\def\coh{\mathop{\text{\upshape{coh}}}} 
\def\coker{\operatorname {coker}} 

\def\Ext{\operatorname {Ext}}

\def\gkdim{\operatorname {GK dim}}

\def\grmod{\operatorname {grmod}} 
 
\def\Hom{\operatorname {Hom}}

\def\im{\operatorname {im}} 

\def\K0{{K_{0}(\blb P^{2}_{q})}} 
\def\ker{\operatorname {ker}}

\def\la{\operatorname{\lambda}} 

\def\P2q{\operatorname {{\blb P}^{2}_{q}}} 

\def\pd{\operatorname {pd}} 

\def\Pic{\operatorname {Pic}} 
\def\PP{\operatorname {\blb P}} 
\def\pr{\mathop{\text{pr}}\nolimits} 

\def\r{\rightarrow} 
 
\def\rank{\operatorname {rank}} 

\def\sgn{\operatorname{sgn}} 
 
\def\Supp{\mathop{\text{\upshape Supp}}} 
 
\def\Tor{\operatorname {Tor}} 

\DeclareMathOperator{\Aut}{Aut}

\DeclareMathOperator{\Div}{Div}

\newtheorem*{theoremA}{Theorem A} 
\newtheorem*{theoremB}{Theorem B}

\newtheorem{lemma}{Lemma}[section] 
\newtheorem{proposition}[lemma]{Proposition} 
\newtheorem{theorem}[lemma]{Theorem} 
\newtheorem{corollary}[lemma]{Corollary}

\theoremstyle{definition} 

\newtheorem{example}[lemma]{Example}

{ 
 
\newtheorem{step}{Step} 
\newtheorem{case}{Case}

} 

\theoremstyle{remark} 

\newtheorem{remark}[lemma]{Remark}

\newdimen\uboxsep \uboxsep=1ex 
\def\uboxn#1{\vtop to 0pt{\hrule height 0pt depth 0pt\vskip\uboxsep 
\hbox to 0pt{\hss #1\hss}\vss}} 

\def\uboxs#1{\vbox to 0pt{\vss\hbox to 0pt{\hss #1\hss} 
\vskip\uboxsep\hrule height 0pt depth 0pt}}

\numberwithin{equation}{section}


\author{Koen De Naeghel} 

\address{Koen De Naeghel \\Departement WNI\\Hasselt University\\ Agoralaan gebouw D \\ B-3590 
Diepenbeek (Belgium).}
\email[K. De Naeghel]{koen.denaeghel@uhasselt.be}

\thanks{The author was supported by LIEGRITS, a Marie Curie Research Training Network funded by the European community as project MRTN-CT 2003-505078}

\date{June 6, 2006} 
\title[Hilbert series of modules of GK-dimension two over elliptic algebras]{
Hilbert series of modules of GK-dimension two over elliptic algebras} 
\begin{document} 

\begin{abstract} 
We characterize the Hilbert functions and minimal resolutions of (critical) Cohen-Macaulay graded right modules of Gelfand-Kirillov dimension two over generic quadratic and cubic three dimensional Artin-Schelter regular algebras.
\end{abstract} 

\maketitle 

\tableofcontents 

\section{Introduction and main results} 

\subsection{Introduction}
In this paper we completely characterize the Hilbert series of critical graded Cohen-Macaulay modules of GK-dimension two over generic elliptic three dimensional Artin-Schelter 
regular algebras which are generated in degree one (see Theorem B 
below).  Such modules can be viewed as irreducible curves in non-commutative deformations of $\PP^2$. 

Our results complete a project started by Ajitabh \cite{Ajitabh,Ajitabh2}. They form a natural counterpart to \cite{DV2,DM} where one describes the possible Hilbert series for modules of GK-dimension three. 

\subsection{Elliptic algebras} \label{Elliptic algebras}

Let $k$ denote the field of complex numbers $\CC$.

We will be dealing with three dimensional Artin-Schelter regular $k$-algebras \cite{AS}. These graded connected algebras were classified in \cite{ATV1,ATV2,Steph1,Steph2} and have all expected nice homological properties. For example they are both left and right noetherian domains with global dimension three and Gelfand-Kirillov dimension three.  

We will furthermore assume $A$ to be generated in degree zero, and we require that $A$ is {\em generic} by which we mean that in the triple $(E,\Lscr,\sigma)$ associated to $A$ \cite{ATV2}, $E$ is a smooth elliptic curve and $\sigma$ is a translation under the group law of $E$ of infinite order. This is equivalent \cite{ATV1} with saying that $A$ takes one of the following forms:
\begin{itemize}
\item
$A$ is quadratic:
\[
A = k \langle x,y,z \rangle/(f_{1},f_{2},f_{3})
\]
where $f_{1},f_{2},f_{3}$ are the homogeneous quadratic relations
\begin{equation*} 
\left\{ 
\begin{array}{l} 
f_{1} = ayz + bzy + cx^{2}  \\ 
f_{2} = azx + bxz + cy^{2}\\ 
f_{3} = axy + byx + cz^{2} 
\end{array} \right. 
\end{equation*} 
where $(a,b,c) \in \PP^{2}$ for which $abc \neq 0$ and $(3abc)^{3} \neq 
(a^{3} + b^{3} + c^{3})^{3}$.
\item
$A$ is cubic:
\[
k\langle x,y \rangle/(g_{1},g_{2})
\]
where $g_{1},g_{2}$ are the homogeneous cubic relations
\begin{equation*} 
\left\{ 
\begin{array}{l} 
g_{1} = ay^{2}x + byxy + axy^{2} + cx^{3}  \\ 
g_{2} = ax^{2}y + bxyx + ayx^{2} + cy^{3}
\end{array} \right. 
\end{equation*} 
where $(a,b,c) \in \PP^{2}$ for which $abc \neq 0$, $b^2 \neq c^2$ and $(2bc)^2 \neq (4a^2 - b^2 - c^2)^2$.
\end{itemize}
In this case $A$ contains a central element $g$ (which is unique up to scalar multiplication), homogeneous of degree three if $A$ is quadratic and of degree four if $A$ is cubic \cite{ATV1}. For later use we put $r_A$ equal to the number of generators of $A$ i.e.
\begin{eqnarray*}
r_A = 
\left\{
\begin{array}{ll}
3 & \text{ if $A$ is quadratic} \\
2 & \text{ if $A$ is cubic}
\end{array}
\right.
\end{eqnarray*}
For the rest of this paper, we will assume $A$ to be such a generic three dimensional Artin-Schelter regular algebra, either quadratic or cubic. In the quadratic case these algebras are so-called three-dimensional Sklyanin algebras (for which the translation $\sigma$ has infinite order). 

\subsection{Some terminology}

By an $A$-module we will mean a finitely generated graded right $A$-module. We write $\grmod(A)$ for the category of $A$-modules. For an $A$-module $M$ and $n \in \ZZ$, write $M_{\leq n} = \bigoplus_{d \leq n}M_{d}$. Define $M(n)$ as the $A$-module equal to $M$ with its original $A$ action, but which is graded by $M(n)_{i} = M_{n + i}$. We refer to the modules $M(n)$ as {\em shifts} ({\em of grading}) of $M$. We say $M$ is {\em normalized} if $M_{< 0} = 0$ and $M_0 \neq 0$. 
The Hilbert series of $M$ is denoted by
\[
h_{M}(t) = \sum_{i = - \infty}^{+ \infty}{(\dim_{k}M_{i})t^{i}} \in\ZZ((t)) 
\]
which makes sense since $A$ is right noetherian. The Hilbert series of $A$ is \cite{AS}
\begin{eqnarray*}
h_A(t) = 
\left\{
\begin{array}{ll}
\frac{1}{(1-t)^3} & \text{ if $A$ is quadratic} \\
\frac{1}{(1-t)^2(1-t^2)} & \text{ if $A$ is cubic}
\end{array}
\right.
\end{eqnarray*}
Taking Hilbert series of a projective (hence free) resolution 
it is easy to see that there exist integers $r,a$ and a Laurent polynomial $s(t) \in \ZZ[t,t^{-1}]$ such that the Hilbert series of $M$ is of the form
\begin{equation} \label{generalHilbertseries}
h_M(t) = h_A(t)(r + a(1-t) - s(t)(1-t)^2)
\end{equation}
We write $\gkdim M$ for the Gelfand-Kirillov dimension (GK-dimension for short) of $M$. As $\gkdim A = 3$, $\gkdim M \leq 3$ and it may be computed as the order of the pole of $h_M(t)$ at $t= 1$, see \cite{ATV2}. The leading coefficient $e_M$ of the series expansion of $h_M(t)$ in powers of $1-t$ is called the {\em multiplicity} of $M$. It is positive and by \eqref{generalHilbertseries} an integer multiple of the multiplicity $e_A$ of $A$, thus as in \cite{ATV2} it will be convenient to put $\iota_A = e_A^{-1}$ and $\epsilon_M = \iota_A e_M$ i.e.
\begin{eqnarray*}
\iota_A = 4 - r_A = 
\left\{
\begin{array}{ll}
1 & \text{ if $A$ is quadratic} \\
2 & \text{ if $A$ is cubic}
\end{array}
\right.
\text{ and }
\epsilon_M = 
\left\{
\begin{array}{ll}
e_M & \text{ if $A$ is quadratic} \\
2e_M & \text{ if $A$ is cubic}
\end{array}
\right.
\end{eqnarray*}
An  $A$-module $M$ is called {\em pure} if for all non-trivial submodules $N \subset M$ we have $\gkdim N = \gkdim M$. If in addition $e_N = e_M$ for all non-trivial submodules we say that $M$ is {\em critical}. This is equivalent with saying that every proper quotient of $M$ has lower GK-dimension. Any pure module $M$ of GK-dimension $d$ admits a filtration such that the successive quotiens are critical of GK-dimension $d$. The graded $\Hom$ and $\Ext$ groups in $\grmod(A)$ will be written as $\underline{\Hom}$ and $\underline{\Ext}$.
We say that $M$ is {\em Cohen-Macaulay} if $\pd M = 3 - \gkdim M$, or equivalently if $\underline{\Ext}^i_A(M,A) = 0$ for $i \neq 3 - \gkdim M$ .

\subsection{Modules of projective dimension one} \label{Modules of projective dimension one}

In this paper we will be concerned with $A$-modules of projective dimension one.
Such a module $M$ admits a minimal resolution of the form
\begin{equation} \label{minres}
0 \r \bigoplus_i A(-i)^{b_i} \rightarrow \bigoplus_i A(-i)^{a_i} \r M \r 0
\end{equation}
The finitely supported sequences of non-negative integers $(a_i)$,$(b_i)$ are usually called the ({\em graded}) {\em Betti numbers} of $M$. Taking Hilbert series of \eqref{minres} one sees they are related to the Hilbert series of $M$ by the formula
\begin{equation} \label{characteristic}
h_M(t) = h_A(t)\sum_i (a_i - b_i)t^i
\end{equation}
The polynomial $q_M(t) = \sum_i (a_i - b_i)t^i \in \ZZ[t]$ is the so-called {\em characteristic polynomial} of $M$. We also write $p_M(t) = q_M(t)/(1-t) \in \ZZ[t,t^{-1}]$.

\subsection{Main results} \label{Modules of GK-dimension two and main results}

For an $A$-module $M$ of GK-dimension two the following assertions are equivalent \cite[\S 4]{ATV2}
\begin{enumerate}
\item
$M$ is pure of projective dimension one,
\item
$M$ has projective dimension one,
\item
$M$ is Cohen-Macaulay,
\item
$M = M^{\vee \vee}$, where $M^{\vee} = \underline{\Ext}^{1}_{A}(M,A)$.
\end{enumerate}
Hence any $A$-module $M$ of GK-dimension two is (uniquely) represented by a pure module of GK-dimension two and projective dimension one, namely $M^{\vee \vee}$.

In order to state our main results, we will first need some terminology \cite[\S 4.1]{DV2}. For positive integers $m,n$ consider the rectangle 
\[ 
R_{m,n} = [1,m] \times [1,n] = 
\{ (\alpha, \beta) \mid 1 \leq \alpha \leq m, 1 \leq \beta \leq n \} 
\subset \ZZ^{2} 
\] 
A subset $L \subset R_{m,n}$ is called a 
{\em ladder} if 
\begin{equation*} 
\forall (\alpha,\beta) \in R_{m,n}: (\alpha,\beta) \not\in L 
\Rightarrow (\alpha+1,\beta), (\alpha,\beta-1) \not\in L 
\end{equation*} 
\begin{example} 
The ladder below is indicated with a dotted line. 

\begin{picture}(0,120) 
\put(40,0){\line(0,5){110}} 
\put(30,100){\line(5,0){100}} 
\put(36.8,0){$\vee$} 
\put(124,97.5){$>$} 
\put(27,5){$\alpha$} 
\put(137,100){$\beta$} 

\qbezier[80](45,95)(85,95)(125,95) 
\qbezier[30](45,95)(45,80)(45,65) 
\qbezier[30](45,65)(60,65)(75,65) 
\qbezier[10](75,55)(75,60)(75,65) 
\qbezier[10](75,55)(80,55)(85,55) 
\qbezier[20](85,35)(85,45)(85,55) 
\qbezier[10](85,35)(90,35)(95,35) 
\qbezier[30](95,5)(95,20)(95,35) 
\qbezier[30](95,5)(110,5)(125,5) 
\qbezier[90](125,5)(125,50)(125,95) 

\put(50,90){$\cdot$} \put(60,90){$\cdot$} \put(70,90){$\cdot$} 
\put(80,90){$\cdot$} \put(90,90){$\cdot$} \put(100,90){$\cdot$} 
\put(110,90){$\cdot$} \put(120,90){$\cdot$} 

\put(50,80){$\cdot$} \put(60,80){$\cdot$} \put(70,80){$\cdot$} 
\put(80,80){$\cdot$} \put(90,80){$\cdot$} \put(100,80){$\cdot$} 
\put(110,80){$\cdot$} \put(120,80){$\cdot$} 

\put(50,70){$\cdot$} \put(60,70){$\cdot$} \put(70,70){$\cdot$} 
\put(80,70){$\cdot$} \put(90,70){$\cdot$} \put(100,70){$\cdot$} 
\put(110,70){$\cdot$} \put(120,70){$\cdot$} 

\put(80,60){$\cdot$} \put(90,60){$\cdot$} \put(100,60){$\cdot$} 
\put(110,60){$\cdot$} \put(120,60){$\cdot$} 

\put(90,50){$\cdot$} \put(100,50){$\cdot$} \put(110,50){$\cdot$} 
\put(120,50){$\cdot$} 

\put(90,40){$\cdot$} \put(100,40){$\cdot$} \put(110,40){$\cdot$} 
\put(120,40){$\cdot$} 

\put(100,30){$\cdot$} \put(110,30){$\cdot$} \put(120,30){$\cdot$} 

\put(100,20){$\cdot$} \put(110,20){$\cdot$} \put(120,20){$\cdot$} 

\put(100,10){$\cdot$} \put(110,10){$\cdot$} \put(120,10){$\cdot$} 
\end{picture} 
\end{example} 
Let $(c_{i})$ be a finitely supported sequence of non-negative integers. We associate a sequence $S(c)$ of length $\sum_i c_i$ to $(c_{i})$ as follows 
\[ 
\ldots\,,\,\underbrace{i-1,\ldots,i-1}_{c_{i-1}\text{ times }}\,,\, 
\underbrace{i,\ldots,i}_{c_i\text{ times 
}}\,,\,\underbrace{i+1,\ldots,i+1}_{c_{i+1}\text{ times }}\,,\,\ldots 
\] 
where by convention the left most non-zero entry of $S(c)$ has index 
one. 

To finitely supported sequences of integers $(a_i)$, $(b_i)$ we associate the matrix $S = S(a,b) = (S(b)_{\beta} - S(a)_{\alpha})_{\alpha \beta}$. It has the properties of a ``degree matrix'':
\begin{equation} \label{degreematrix}
S_{\alpha+1,\beta} \leq S_{\alpha \beta} \leq S_{\alpha,\beta+1}\quad \text{ and } \quad 
S_{\alpha \beta} - S_{\alpha \beta'} = S_{\alpha' \beta} - S_{\alpha' \beta'}
\end{equation}
from which it follows that  
\begin{equation} \label{ref-5.1-30} 
L_{a,b}= \{(\alpha,\beta)\in R_{m,n} \mid S(a)_\alpha< S(b)_\beta\} 
\text{ where } m = \sum_i a_i, n = \sum_i b_i 
\end{equation} 
is a ladder. Our following main result is proved in \S\ref{Proof of Theorem A} below.
\begin{theoremA} \label{Betti}
Let $(a_i)$, $(b_i)$ be finitely supported sequences of integers and put $m = \sum_i a_i$, $n = \sum_i b_i$. 
\begin{enumerate}
\item
$(a_i)$, $(b_i)$ appear as the Betti numbers of a 
graded right $A$-module $M$ of GK-dimension two and projective dimension one if and only if 
\begin{enumerate} 
\item The $(a_i),(b_i)$ are non-negative. 
\item $\sum_i b_i = \sum_i a_i$. 
\item $\forall (\alpha, \beta) \in  R_{n,m}: \beta \geq \alpha \Rightarrow (\alpha, \beta)  \in L_{a,b}$. 
\end{enumerate} 
\item
$(a_i)$, $(b_i)$ appear as the Betti numbers of a \underline{critical} graded right $A$-module $M$ of GK-dimension two and projective dimension one if and only if 
\begin{enumerate} 
\item The $(a_i),(b_i)$ are non-negative. 
\item $\sum_i b_i = \sum_i a_i$. 
\item $\forall (\alpha, \beta) \in  R_{n,m}: \beta \geq \alpha - 1 \Rightarrow (\alpha, \beta)  \in L_{a,b}$. 
\item
If $A$ is cubic it is not true that {\em (}$n\geq 2$ and $\forall \alpha, \beta: S(b)_{\beta}- S(a)_{\alpha} = 1${\em )}.
\end{enumerate} 
\end{enumerate}
In both statements, the module $M$ may be chosen to be $g$-torsion free.
\end{theoremA}
\begin{remark}
For quadratic $A$ it was proved in \cite{Ajitabh} that the appearing conditions in Theorem A are necessary, and in \cite{Ajitabh2} there were shown to be sufficient in the case where $\sum_i b_i = \sum_i a_i = 1$.  
\end{remark}
\begin{remark}
Theorem A is an analogue of the description of the Betti numbers of pure $A$-modules of GK-dimension three and projective dimension one, see \cite{DV2,DM}.
\end{remark}
By Theorem A, a minimal resolution of a (resp. critical) $A$-module $M$ of GK-dimension two and projective dimension one is of the form 
\[
0 \r \bigoplus_i A(-i)^{b_i} \r \bigoplus_i A(-i)^{a_i} \r M \r 0
\]
for which a generic map $f: \bigoplus_i A(-i)^{b_i} \r \bigoplus_i A(-i)^{a_i}$ is represented by left matrix multiplication with a matrix of the form 
\begin{eqnarray*} 
f = 
\begin{pmatrix}
\ast  & \ast & \ast & \dots  & \ast \\
& \ast & \ast & \dots & \ast  \\
 & & \ast & \dots  & \ast  \\
 & &  & \ddots & \vdots \\
 &  &  & & \ast 
\end{pmatrix}
\text{ resp. }
f = 
\begin{pmatrix}
\ast & \ast  & \ast  & \dots  & \ast \\
\ast & \ast & \ast & \dots & \ast  \\
 & \ast & \ast & \dots  & \ast  \\
 & &  \ddots  & \ddots & \vdots \\
 &  &  & \ast & \ast 
\end{pmatrix}
\end{eqnarray*}
where the indicated entries $\ast$ are nonzero homogeneous elements in $A$ of degree $\geq 1$. In case $A$ is cubic, there is an additional condition for critical $M$: Not all entries in $f$ have degree $1$, unless $f$ is a $1 \times 1$ matrix (reflecting condition (2)(d) in Theorem A). In other words, in case $A$ is cubic then the minimal resolution of a critical graded right $A$-module $M$ of GK-dimension two and projective dimension one cannot be the form (up to shift of grading)
\[
0 \r A(-1)^n \r A^n \r M \r 0, \quad n \geq 2
\]
This might seem surprising for the reader. The reason is explained in Example \ref{examplestriking} below.

As a concequence of Theorem A, we will deduce in \S\ref{Proof of Theorem B and other properties of Hilbert series}
\begin{theoremB} \label{Hilbert}
Let $\epsilon > 0$ be an integer and put $e = \epsilon/\iota_A$. There is a bijective correspondence between Hilbert series $h(t)$ of normalized $A$-modules of GK-dimension two, projective dimension one and multiplicity $e$, and polynomials $s(t) = \sum_{i}s_{i}t^{i} \in \ZZ[t]$ which satisfy
\begin{equation} \label{Hilbert1}
\epsilon > s_{0} \geq s_{1} \geq \dots \geq 0
\end{equation}
The correspondence is given by $h(t) = h_A(t)(\epsilon(1-t)-s(t)(1-t)^2)$, explicitely
\begin{eqnarray} \label{Hilbertseries}
h(t) 
=
\left\{
\begin{array}{ll}
\frac{e}{(1-t)^{2}} - \frac{s(t)}{1-t} & \text{ if $A$ is quadratic} \\
\frac{2e}{(1-t)(1-t^{2})} - \frac{s(t)}{1-t^2} & \text{ if $A$ is cubic} \\
\end{array}
\right.
\end{eqnarray}
Further, if we restrict to critical $A$-modules then the same statement holds where \eqref{Hilbert1} is replaced by
\begin{equation} \label{Hilbert2}
\epsilon > s_{0} > s_{1} > \dots \geq 0
\text{ and if $A$ is cubic and $\epsilon > 1$ then $s(t) \neq 0$ }
\end{equation}
\end{theoremB}
It is clear that there are only finitely many polynomials $s(t) \in \ZZ[t]$ which satisfy \eqref{Hilbert2}. Hence Theorem B implies that
there are only finitely many possibilities for the Hilbert series of a critical normalized Cohen-Macaulay $A$-module of GK-dimension two and multiplicity $e$. This consequence was already observed by Ajitabh in \cite{Ajitabh} for quadratic $A$. In fact it is easy to count the number of possibilities.
\begin{corollary} \label{count}
Let $\epsilon > 0$ be an integer and put $e = \epsilon/\iota_A$. The number of Laurent power series which appear as the Hilbert series of a critical normalized module of GK-dimension two, projective dimension one and multiplicity $e$ is equal to 
\begin{eqnarray*}
\left\{
\begin{array}{ll}
2^{\epsilon-1} - 1 & \text{ if $A$ is cubic and $\epsilon > 1$,} \\
2^{\epsilon-1} & \text{ else. }
\end{array}
\right.
\end{eqnarray*}
\end{corollary}

\begin{remark}
It follows from Theorem B that there are infinitely many possibilities for the Hilbert series of a normalized Cohen-Macaulay module of GK-dimension two and multiplicity $e > 1$. This is to be expected, for example if $A$ is quadratic and $S$ is a line module over $A$ then $M = S^{e-1} \oplus S(-n)$ is a (non-critical) normalized Cohen-Macaulay $A$-module GK-dimension two and multiplicity $e$, for all integers $n \geq 0$.
\end{remark}
\begin{remark}
For the convenience of the reader we have included in
Appendix \ref{A} the list of possible Hilbert series and Betti numbers of critical normalized Cohen-Macaulay modules $M$ of GK-dimension two and $\epsilon_M \leq 4$. 
\end{remark}
\begin{remark}
It is well-known that Theorem A (and hence Theorem B) holds for the commutative polynomial ring $k[x,y,z]$, which is a non-generic quadratic three dimensional Artin-Schelter regular algebra. See for example \cite[Proposition 2.7 and Theorem 2.8]{BCG}. We conjecture that Theorems A and B are true for all three dimensional Artin-Schelter regular algebras generated in degree one (thus not only the generic ones). 
\end{remark}
We end this introduction by saying a few words about the proof of Theorem A. The most difficult part is to show that the conditions in Theorem A(2) are sufficient \S\ref{Proof that the conditions in Theorem A(2) are sufficient}. Roughly, this will be derived from the following three observations:
\begin{itemize}
\item
To any $g$-torsionfree $A$-module $M$ of GK-dimension two one may associate a divisor on the elliptic curve $E$, denoted by $\Div M$. This notion was introduced by Ajitabh in \cite{Ajitabh}, who showed in \cite{Ajitabh2} that writing $\Div M = D + (q)$ for some effective divisor $D$ on $E$, there is a sufficient condition on $D$ (called quantum-irreducibility) for $M$ to be critical. See also \S\ref{The divisor of a curve module} and \S\ref{Quantum-irreducible divisors on $E$} below.
\item
For any positive integer, there exists an effective multplicity-free quantum-irreducible divisor on $E$. This was shown in \cite{Ajitabh2}. See also \S\ref{Quantum-irreducible divisors on $E$} below.
\item
Let $(a_i)$, $(b_i)$ be finitely supported sequences of integers satisfying Theorem A(2). Let $D$ be a multiplicity-free effective divisor on $E$ of degree $r_A \sum_i i(b_i - a_i) - 1$. In Theorem \ref{strongthm} below we show that there is a $g$-torsion free $M \in \grmod(A)$ of GK-dimension two which has a minimal resolution 
\[
0 \r \bigoplus_i A(-i)^{b_i} \xrightarrow{f} \bigoplus_i A(-i)^{a_i} \r M \r 0,
\]
for which the matrix represening the map $f$ has the form  
\begin{eqnarray*} 
f = 
\begin{pmatrix}
\ast & \ast  & \ast & \dots & \ast & \ast \\
\ast &  & &  & & \ast  \\
 & \ast & &  & & \ast  \\
 & &  \ddots  &  & & \vdots \\
 &  &  & \ast &  & \ast \\
 &  &  & & \ast & \ast 
\end{pmatrix}
\end{eqnarray*}
(where the entries off the diagonal, first row and last column are zero) and for which $\Div M = D + (q)$ for some $q\in E$.  
\end{itemize}

\subsection{Acknowledgements}

I would like to state my gratitude to Michel Van den Bergh for his useful coments on a preliminary version of this paper. I also thank Elisa Gorla for explaining the analogy of our results in the commutative case. 

\section{Preliminaries} \label{Preliminaries}

Throughout we will assume $A$ to be a generic three-dimensional Artin-Schelter regular algebra, either quadratic or cubic, as decribed in \S\ref{Elliptic algebras}. 

\subsection{Geometric data}

In this part we recall some terminology and basic facts on elliptic algebras from \cite{ATV1,ATV2}. 

The algebra $A$ is completely determined by geometric data $(E,\Lscr,\sigma)$ where 
\begin{itemize}
\item
if $A$ is quadratic then $j: E \hookrightarrow \PP^2$ is a divisor of degree three, $\Lscr = j^{\ast}\Oscr_{\PP^2}(1)$ line bundle of degree three and $\sigma \in \Aut(E)$, 
\item
if $A$ is cubic then $j: E \hookrightarrow \PP^1 \times \PP^1$ is a divisor of bidegree $(2,2)$, $\Lscr = j^{\ast}\pr_1^{\ast}\Oscr_{\PP^1}(1)$ line bundle of degree two and $\sigma \in \Aut(E)$.
\end{itemize}
As we choose $A$ to be generic, $E$ is smooth curve of arithmetic genus one i.e. an elliptic curve, and $\sigma$ is a translation on $E$. In case $A$ is cubic then $\sigma$ is of the form $\sigma(q_1,q_2) = (q_2,f(q_1,q_2))$ for some map $f: E \r \PP^1$. 

Let $\Escr \in \Pic(E)$ be a line bundle on $E$. We use the notation $\Escr^{\sigma}$ for the pull-back $\sigma^{\ast}\Escr$. 
Thus $(\Escr^\sigma)_p = \Escr_{p^{\sigma}}$ for $p \in E$.
We regard $\Pic(E)$ as a module over $\ZZ[\sigma,\sigma^{-1}]$, where the action of a Laurent polynomial $f(\sigma) = \sum_{i}a_i \sigma^i$ on $\Escr$ is defined as 
\[
\Escr^{f(\sigma)} := \otimes_i (\Escr^{\sigma^{i}})^{\otimes a_i}
\]
Recall that there is, up to scalar multiplication, a canonical central element $g \in A$, homogeneous of degree $\iota_A r_A$. The factor ring ring $A/gA$ is isomorphic to the twisted homogeneous coordinate ring 
\[
B = \bigoplus_{n \geq 0}H^{0}(E,\Lscr_n) \text{ where } 
\Lscr_n = \Lscr \otimes_E \Lscr^{\sigma} \otimes_E \dots \otimes_E \Lscr^{n-1} = \Lscr^{(1-\sigma^n)/(1-\sigma)}
\]
is a line bundle of degree $r_A n$. Note that in case $A$ is cubic we have $\Lscr = j^\ast\pr_2^{\ast}\Oscr_{\PP^1}(1)$. Multiplication in $B$ is defined by $b_n b_m = b_{n} \otimes_E b_{m}^{\sigma^{n}}$ for $b_n \in B_n$, $b_m \in B_m$, where $b_m^{\sigma^{n}} = b_m \circ \sigma^n$. The algebra $B$ has Gelfand-Kirillov dimension two, and it is a domain since $E$ is reduced. The homogeneous elements of $B$ will be identified with the corresponding sections of the appropriate line bundles on $E$. For any $m \in A_n$, we denote by $\overline{m}$ its image in $B \cong A/Ag$.

There is a (left exact) global section functor
\[
\Gamma_{\ast}: \coh(E) \r \grmod(B): \Fscr \r \bigoplus_{n \geq 0} H^{0}(E,\Fscr \otimes_E \Lscr_n)
\]
whose right adjoint is exact, denoted by $\widetilde{(-)}$.
It was shown in \cite{AVdB} that they induce a category equivalence between $\coh(E)$ and $\grmod(B)/\grmod(B)_0$. Here $\grmod(B)_0$ stands for the Serre subcategory of the finite length modules in the category $\grmod(B)$ of finitely generated graded right $B$-modules. 

It will be convenient below to let the shift functors $-(n)$ on $\coh(E)$ be the ones obtained from the equivalence and not the ones coming from the embedding $j$. Thus $\Oscr_E(n) = \sigma_{\ast}^n \Lscr_n$ and $\Fscr(n) = \sigma_{\ast}^n\Fscr \otimes_E \Oscr_E(n)$ for $\Fscr \in \coh(E)$. 

For $p \in E$ we write $P = \left(\Gamma_{\ast}(k(p))\right)_A \in \grmod(A)$ where $k(p)$ is the skyscraper sheaf $k$ sitting at $p$. Observe that $k(p)(n) = k(p^{\sigma^n})$. Such $A$-modules $P$ are called {\em point modules} over $A$. It is easy to see that $h_{P}(t) = (1-t)^{-1}$. 

\subsection{Group law and divisors on $E$} \label{Group law} 

Fixing a group law on $E$ the automorphism $\sigma$ is a translation by some point $\xi \in E$. Thus $p^{\sigma} = p + \xi$ for $p \in E$. We write $o$ for the orgin of the group law. Linear equivalence of divisors $D,D'$ on $E$ will be denoted by $D \sim D'$. We will frequently use 
\begin{proposition} {\em(}\cite[IV Theorem 4.13B]{H}{\em)} \label{H}
Let $D,D'$ be two divisors on $E$. Then
\[
D \sim D' \Leftrightarrow \deg D = \deg D' \text{ and $D,D'$ have the same sum in the group law of $E$}   
\]
\end{proposition}
For example, for three points $p,q,r \in E$ we have $p = q + r$ in the group law of $E$ if and only if $(p) + (o) \sim (q) + (r)$ as divisors on $E$.

For a nonzero global section $s \in H^{0}(E,\Escr)$ and $p \in E$ we write 
$s(p)$ for the image of $s$ in the one dimensional $k$-linear space $\Escr \otimes_E k(p) \cong \Escr_p / m_p \Escr_p$, where $m_p$ is the maximal ideal of the local ring $\Oscr_p$ and $k(p) = \Oscr_p / m_p$. In case $\Escr = \Oscr_E(i)$ we have $s(p) \in \Oscr_E(i) \otimes_E k(p) = k(p^{\sigma^{i}})$.

We write $s^{\sigma}$ for the image of $s$ under the $k$-linear isomorphism $H^{0}(E,\Escr) \cong H^{0}(E,\Escr^{\sigma})$. We have $s^{\sigma}(p) = s(p^\sigma)$ under the isomorphism $k(p)^{\sigma} \cong k(p^{\sigma})$. We write $\Div(s)$ for the divisor of zeros of $s$. It follows that $\Div(s^{\sigma}) = \sigma^{-1}\Div(s)$. 

Consider a map $N:\Oscr_E(-j) \r \Oscr_E(-i)$ where $i < j$. As $\Oscr_E(-i) \otimes k(p^{\sigma^{i}}) = k(p)$, a point $p \in E$ is supported on the cokernel of $N = (n)$ if and only if $N \otimes k(p^{\sigma^{i}}) =  (n^{\sigma^{i}}(p))$ is zero. Here, $n$ is viewed as a global section of $\Oscr_E(j-i)$.

This is generalized as follows. Let $(a_i)$, $(b_i)$ be finitely supported sequences of non-negative integers. Consider a map $N: \bigoplus_i \Oscr_E(-i)^{b_i} \r \bigoplus_i \Oscr_E(-i)^{a_i}$. To $N = (n_{\alpha \beta})_{\alpha \beta}$ we associate a new matrix $X_N$, given by $(X_{N})_{\alpha \beta} = n^{{\sigma^{S(a)_\alpha}}}_{\alpha\beta}$. It is easy to see that a point $p \in E$ is supported on the cokernel of $N = (n_{\alpha \beta})_{\alpha \beta}$ if and only if the rank of the matrix $X_N(p)$ is less than $\sum_i a_i$, where 
\[
X_N(p) := X_N \otimes k(p) \quad \text{ i.e. } \left(X_N(p)\right)_{\alpha \beta} = 
n^{{\sigma^{S(a)_\alpha}}}_{\alpha\beta}(p) 
\]

\subsection{The divisor of a curve module} \label{The divisor of a curve module}

By a {\em curve $A$-module} \cite{Ajitabh} we will mean a $g$-torsion free module $A$-module $M$ of GK-dimension two. 

It was shown in \cite{Ajitabh, AjVdB} that to any curve $A$-module one may associate a divisor on $E$. Actually this was done in case $A$ is quadratic, but a similar treatment holds for cubic $A$. Let us recall how this is done, considering both cases (quadratic and cubic) at the same time.

Let $M$ be a curve $A$-module. As $M$ is $g$-torsion free, $M/Mg$ has GK-dimension one. Hence $(M/Mg){\,\widetilde{}\,\,}$ is a finite dimensional $\Oscr_E$-module which corresponds to a divisor on $E$. We will call this the {\em divisor} of $M$ and denote it by $\Div(M)$. 
\begin{proposition} {\em(}\cite{AjVdB}{\em)} \label{propsDiv}
Let $M$ be curve $A$-module.
\begin{enumerate}
\item
$\Div(M)$ is an effective divisor of degree $r_A\epsilon_M$.
\item
For any integer $l$ we have $\Div(M(l)) = \sigma^l \Div(M)$.
\item
$\Div$ is additive on short exact sequences i.e. for a short exact sequence of curve modules $0 \r M' \r M \r M'' \r 0$ in $\grmod(A)$ we have
\[
\Div(M) = \Div(M') + \Div(M'')
\]
\item
Let $p \in E$ and write $P = \left(\Gamma_{\ast}(\Oscr_p)\right)_A$ for the corresponding point module. Assume we have an exact sequence $0 \r K \r M \xrightarrow{f} P$ where $f \neq 0$. Then $K$ is a curve $A$-module and 
\[
\Div(K) = \Div(M) - (p) + (p^{\sigma^{-\iota_A r_A}})
\]
\item
Let $p \in E$. If $\Hom_A(M,P) \neq 0$ then $p \in \Supp(\Div(M))$. In case $M_{<0} = 0$ the converse is also holds.
\end{enumerate}
\end{proposition}
We also mention
\begin{lemma} \label{purecohmac}
Let $M$ be a pure curve $A$-module. Then $\Div(M) = \Div(M^{\vee \vee})$.
\end{lemma}
\begin{proof}
By \cite[Corollary 4.2]{ATV2} the canonical map $\mu_M: M \r M^{\vee \vee}$ is injective and its cokernel is finite dimensional. Thus $\pi M = \pi M^{\vee\vee}$ and hence $(M/Mg){\,\widetilde{}\,\,} = (M^{\vee\vee}/M^{\vee\vee}g){\,\widetilde{}\,\,}$. This means that $\Div(M) = \Div(M^{\vee \vee})$.
\end{proof}
For any $g$-torsion free $a \in A_n$ the divisor $\Div(A/aA)$ of the curve $A$-module $M = A/aA$ coincides with the divisor of zeros $\Div(\overline{a})$ of the global section $\overline{a} \in H^{0}(E,\Lscr_n)$. Indeed, this follows from the short exact sequence in $\coh(E)$
\[
0 \r \Oscr_E(-n) \xrightarrow{\overline{a}} \Oscr_E \r (M/Mg){\,\widetilde{}\,\,} \r 0
\]
More generally, in \cite{Ajitabh} it was shown that for any Cohen-Macaulay curve $A$-module $M$ we may interpret $\Div(M)$ as the divisor of zeros of some global section $s_{[M]}$ of the invertible sheaf $\Lscr^{p_M(\sigma)}$ on $E$. As this will play a key role further on, we will now recall the construction of $s_{[M]}$. For more details the reader is referred to \cite{Ajitabh, Ajitabh2}.

Let $M$ be a Cohen-Macaulay curve $A$-module, say with minimal resolution
\begin{equation} \label{resmin}
0 \r \bigoplus_i A(-i)^{b_i} \xrightarrow{f} \bigoplus_i A(-i)^{a_i} \r M \r 0
\end{equation}  
We represent the map $f$ in \eqref{resmin} by left multiplication by a matrix $[M]$ whose entries are homogeneous elements $m_{\alpha\beta}$ in $A$. 
Applying the functor $- \otimes_A B$ to \eqref{resmin} we find an exact sequence in $\grmod(B)$
\begin{equation} \label{resmin'}
0 \r \bigoplus_i B(-i)^{b_i} \xrightarrow{\overline{f}} \bigoplus_i B(-i)^{a_i} \r M/Mg \r 0
\end{equation}
where we have used the $g$-torsionfreeness of $M$ to derive 
\[
\Tor_1^{A}(M,A/Ag) = \ker (M(-\iota_Ar_A) \xrightarrow{\cdot g} M) = 0
\] 
The map $\overline{f}$ is represented by $\overline{[M]}$, the matrix obtained from $[M]$ by replacing the entries $m_{\alpha\beta} \in A$ by $\overline{m}_{\alpha\beta} \in B$. Applying the exact functor $(\widetilde{-})$ on \eqref{resmin'} we obtain an exact sequence in $\coh(E)$
\begin{equation} \label{resmin''}
0 \r \bigoplus_i \Oscr_E(-i)^{b_i} \xrightarrow{\overline{[M]}} \bigoplus_i \Oscr_E(-i)^{a_i} \r (M/Mg){\,\widetilde{}\,\,} \r 0
\end{equation}
It is now clear that the divisor of $(M/Mg){\,\widetilde{}\,\,}$ is precisely the zerodivisor of $\det(X_{[M]})$, where as in \S\ref{Group law} the matrix $X_{[M]}$ is defnined as
\[
X_{[M]} = (\overline{m}^{\sigma^{S(a)_\alpha}}_{\alpha\beta})_{\alpha\beta} \text{ where } \overline{m}^{\sigma^{S(a)_\alpha}}_{\alpha\beta} \in H^0(E,\Lscr_{S(b)_\beta - S(a)_\alpha}^{\sigma^{S(a)_\alpha}})
\]
and
\[
\det(X_{[M]}) = \sum_{\gamma \in S_{r}} \sgn(\gamma)
\overline{m}_{1 \gamma(1)}^{\sigma^{S(a)_1}} \otimes \overline{m}_{2 \gamma(2)}^{\sigma^{S(a)_2}} \otimes \dots \otimes \overline{m}_{r \gamma(r)}^{\sigma^{S(a)_r}} 
\]
We denote $s_{[M]} = \det X_{[M]}$.
\begin{example} \label{example}
Consider a Cohen-Macaulay curve $A$-module with minimal resolution of the form 
\[
0 \r A(-3) \oplus A(-7) \oplus A(-8) \xrightarrow{
[M] \cdot} A(-1) \oplus A(-2) \oplus A(-7) \r M \r 0
\] 
where
\[
[M] = 
\begin{pmatrix}
m_{11} & m_{12} & m_{13} \\
m_{21} & m_{22} & m_{23} \\
0 & 0 & m_{33}
\end{pmatrix}
\]
where the entries $m_{ij} \in A$ are homogeneous elements with appropriate degrees.
The corresponding exact sequence on $\coh(E)$ becomes
\[
0 \r \Oscr_E(-3) \oplus \Oscr_E(-7) \oplus \Oscr_E(-8) \xrightarrow{
\overline{[M]}\cdot} \Oscr_E(-1) \oplus \Oscr_E(-2) \oplus \Oscr_E(-7) \r (M/Mg){\,\widetilde{}\,\,} \r 0
\] 
Hence
\[
X_{[M]} = 
\begin{pmatrix}
\overline{m}_{11}^{\sigma} & \overline{m}_{12}^{\sigma} & \overline{m}_{13}^{\sigma} \\
\overline{m}_{21}^{\sigma^{2}} & \overline{m}_{22}^{\sigma^2} & \overline{m}_{23}^{\sigma^2} \\
0 & 0 & \overline{m}_{33}^{\sigma^{7}}
\end{pmatrix} 
\]
Therefore 
\[
s_{[M]} = \det X_{[M]} = \overline{m}_{11}^{\sigma} \otimes \overline{m}_{22}^{\sigma^2} \otimes \overline{m}_{33}^{\sigma^{7}} - \overline{m}_{12}^{\sigma} \otimes \overline{m}_{21}^{\sigma^2} \otimes \overline{m}_{33}^{\sigma^{7}} 
\]
which is a global section of the line bundle $\Lscr^{p_M(\sigma)}$ where $p_M(t) = q_M(t)/(1-t) = t + 2t^2 + t^3 + t^4 + t^5 + t^6 + t^7$, see \S\ref{Modules of projective dimension one}. 
\end{example}
In general, the following result was shown in \cite{Ajitabh,Ajitabh2}.
\begin{proposition} \label{propsM}
Let $M$ be a Cohen-Macaulay curve $A$-module, say with a minimal resolution 
\begin{equation*} 
0 \r \bigoplus_i A(-i)^{b_i} \xrightarrow{[M]\cdot} \bigoplus_i A(-i)^{a_i} \r M \r 0
\end{equation*}
Then the following holds.
\begin{enumerate}
\item
Up to a scalar multiple, $s_{[M]}$ is nonzero and independent of the choice of a minimal resolution for $M$.
\item
For any integer $l$ we have $s_{[M(l)]} = \sigma^l s_{[M]}$.
\item
$s_{[M]} \in H^{0}(E,\Lscr^{p_M(\sigma)})$ and $\deg \Lscr^{p_M(\sigma)} = r_A \epsilon_M$. 
\item
The divisor of zeros of $s_{[M]}$ coincides with the divisor $\Div(M)$.
\end{enumerate}
\end{proposition}

\subsection{Further properties of divisors of curve modules}

By Proposition \ref{propsM}, the divisor of a Cohen-Macaulay curve $A$-module $M$ is the divisor of a section of the line bundle $\Lscr^{p_M(\sigma)}$, and this line bundle depends only on the Hilbert series of $M$. This yields (see also \cite{Ajitabh2})
\begin{proposition} \label{linequiv}
\begin{enumerate}
\item
Let $M,M'$ be two Cohen-Macaulay curve modules with the same Hilbert series $h_{M}(t) = h_{M'}(t)$. Then $\Div(M) \sim \Div(M')$.
\item
Let $D$ be a divisor on $E$. Then, for any Laurent power series $h(t) \in \ZZ((t))$ there is at most one $q \in E$ such that $D + (q) = \Div(M)$ for some Cohen-Macaulay curve module $M$ with Hilbert series $h(t)$.
\end{enumerate}
\end{proposition}
\begin{proof}
(1) 
As $h_M(t) = h_{M'}(t)$ we also have $p_M(t) = p_{M'}(t)$. Proposition \ref{propsM} implies that $\Div(M)$ and $\Div(M')$ are both divisors of global sections of the same line bundle $\Lscr^{p_M(\sigma)}$. 
Hence $\Div(M)$ and $\Div(M')$ are lineary equivalent. 


(2) For the second statement, assume $q,q' \in E$ for which 
\[
D + (q) = \Div(M), \quad D + (q') = \Div(M') 
\] 
for some Cohen-Macaulay curve modules $M,M'$ with Hilbert series $h(t)$. By the first part of the current proposition and Proposition \ref{H}, $\Div(M)$ and $\Div(M')$ have the same sum in the group law of $E$. But this implies $q = q'$, ending the proof.
\end{proof}
\begin{remark}
In Theorem \ref{strongthm} we prove a converse of Proposition \ref{linequiv}(2). As mentioned in the introduction, 
this will be our key result to prove Theorem A.
\end{remark}
In case of Cohen-Macaulay curve $A$-modules of the form $M = A/aA$ we have a more detailed version.
\begin{lemma} \label{dimdiv}
Let $\epsilon > 0$ be an integer and $D$ a multiplicity-free effective divisor of degree $\leq r_A \epsilon$. Then 
\[
\Dscr = \{ b \in B_n \mid \Supp(D) \subset \Supp(\Div(b))
\} 
\]
is a $k$-linear subspace of $B_n$ of dimension 
\begin{eqnarray*}
\dim_k \Dscr 
\left\{
\begin{array}{ll}
= r_A n - \deg D & \text{ if } \deg D \leq r_A n-1 \\
\leq 1 & \text{ if } \deg D = r_A n
\end{array}
\right.
\end{eqnarray*}
\end{lemma}
\begin{proof}
This follows from the category equivalence $\Gamma_\ast$ and Riemann-Roch on the elliptic curve $E$.
\end{proof}
From the previous lemma it is clear that, given $r_A n-2$ points on $E$, there are infinitely many sections on $\Lscr_n$ vanishing in these points. We will need a somewhat more refined version of this.
\begin{lemma} \label{multfree}
Let $n > 0$ be an integer, $q_1,\dots,q_{r_An-2} \in E$ be different points and $\Qscr$ a finite set of points of $E$. There exists a homogeneous form $b \in B_n$ for which $\Div(b)$ is multiplicity-free, $q_i \in \Supp \Div(b)$ and $\Supp \Div(b) \cap \Qscr = \emptyset$.
\end{lemma}

\subsection{Division in $B$}

The following lemma is a useful criterion for division in the twisted homogeneous coordinate ring $B = \Gamma_{\ast}(\Oscr_E)$.
\begin{lemma} \label{division}
Let $b \in B_n$ and $\tilde{b} \in B_m$ be nonzero. Then
\[
\Div(b) = \Div(\tilde{b}) + D \text{ for some effective divisor $D$} \Leftrightarrow b = \tilde{b}c \text{ for some } c \in B_{n-m}.
\]
\end{lemma}
\begin{proof}
Recall that $B_n = H^{0}(E,\Lscr_n)$ where $\Lscr_n$ is the invertible sheaf 
\[
\Lscr_n = \Lscr \otimes \Lscr^{\sigma} \otimes \dots \otimes \Lscr^{\sigma^{n-1}} = \Oscr_E(n)
\]
Write $\Lscr_n = \Oscr_E(D_n)$ for some divisor $D_n$ on $E$. Using this notation, $\Div(b) \sim D_n$ and $\Div(\tilde{b}) \sim D_m$. It follows that $D \sim D_n - D_m$ where 
\[
\Oscr_E(D_n - D_m) \cong \Lscr_n \otimes \Lscr_m^{-1} = \Lscr^{\sigma^m} \otimes \Lscr^{\sigma^{m + 1}} \otimes \dots \otimes \Lscr^{\sigma^{n-1}} = \Lscr_{n-m}^{\sigma^{m}}
\]
As $D$ is effective there is a $c' \in H^{0}(E,\Lscr_{n-m}^{\sigma^{m}})$ for which $\Div(c') = D$. Thus
\[
\Div(b) = \Div(\tilde{b}) + \Div(c') = \Div(\tilde{b} \otimes c') = \Div(\tilde{b}c'^{\sigma^{-m}})
\]
By \cite[II Proposition 7.7]{H} we have $b = \la \tilde{b}c'^{\sigma^{-m}}$ for some $0 \neq \la \in k$. Putting $c = \la c'^{\sigma^{-m}} \in H^0(E,\Lscr_{n-m}) = B_{n-m}$ proves what we want.
\end{proof}
%
%
From the previous lemma we deduce
\begin{lemma} \label{division'}
Let $b, \tilde{b} \in B_n$ be nonzero. Assume 
$\Div(b) = \Div(\tilde{b}) - (p) + (q)$ for some $p,q \in E$. Then $p = q$ and $b = \tilde{b}c$ for some $c \in k$.
\end{lemma} 
\begin{proof}
By Proposition \ref{linequiv}, $\Div(b) \sim  \Div(\tilde{b})$. As these divisors have the same degree, Proposition \ref{H} gives they have the same sum for the group law of $A$. Thus $p = q$. Invoking Lemma \ref{division} completes the proof. 
\end{proof}

\subsection{Quantum-irreducible divisors on $E$} \label{Quantum-irreducible divisors on $E$}

Let $M$ be a curve $A$-module. In \cite{Ajitabh2} the author found a sufficient condition on $\Div(M)$ for $M$ to be critical. We will need this result. For convenience we briefly recall his treatment.

As $\Div$ is additive on short exact sequences (Proposition \ref{propsDiv}) we have
\begin{lemma} {\em(}\cite[Lemma 3.3]{Ajitabh2}{\em)} \label{split2}
Let $M$ be a curve $A$-module. If $M$ is not critical then 
\[
\Div(M) = \Div(M_1) + \Div(M_2)
\]
for some curve $A$-modules $M_1, M_2$.
\end{lemma}
Inspired by the previous lemma, we say that an effective divisor $D$ on $E$ is {\em quantum-reducible} \cite{Ajitabh2} if
\[
D = \Div(M) + D'
\] 
where $M$ is a curve $A$-module and $D'$ is an effective divisor of degree $> 0$. We say $D$ is {\em quantum-irreducible} if $D$ is not quantum-reducible. By Proposition \ref{propsDiv}, any effective divisor of degree $< r_A$ is quantum-irreducible. We have
\begin{lemma} \label{sufficient}
Let $M$ be a curve $A$-module. Assume $\Div(M) = D + (q)$ for some quantum-irreducible divisor $D$ and $q \in E$. Then $M$ is critical.
\end{lemma}
\begin{proof}
Assume by contradiction that $\Div(M)$ is not critical. By Lemma \ref{split2} we have $\Div(M) = \Div(M_1) + \Div(M_2)$ for some curve $A$-modules $M_1, M_2$. Since $\Div(M) = D+(q)$ we must have $q \in \Supp M_i$ for some $i = 1,2$, say for $i = 2$. Then $D' = \Div(M_2) - (q)$ is effective and of degree $> 0$ by Proposition \ref{propsDiv}(1). Now $D = \Div(M_1) + D'$ contradicts the irreducibility of $D$. 
\end{proof}
The existence of quantum-irreducible divisors follows from (it is straightforward to extend the proof for cubic $A$)
\begin{theorem} {\em (}\cite[Theorem 3.7]{Ajitabh2}{\em )} \label{existence}
For any positive integer $n$ there exists a \linebreak  multiplicity-free quantum-irreducible effective divisor $D$ of degree $n$ on $E$, which is not the divisor of a curve $A$-module.
\end{theorem}
\begin{proof}[Sketch of the proof]
It is sufficient to construct a multiplicity-free effective divisor $D$ which is not of the form $D = D' + D''$ for some effective divisors $D',D''$ where $D' \sim \Div M$ for some critical curve module $M \in \grmod(A)$. By Proposition \ref{propsM}, it is sufficient to exclude those $D'$ for which $\Oscr_E(D') \cong \Lscr^{\sigma^n p_M(\sigma)}$ for some critical normalized curve $A$-module $M$ and integer $n$. By \cite{Ajitabh} there are only finitely many possibilities for such $p_M(t) \in \ZZ[t,t^{-1}]$, as there are only finitely many possibilities for the Hilbert series for $M$. This is also part of Theorem B (a part for which we do not rely on the current theorem). Thus we have to exclude a countable number of divisors. As $k = \CC$ is uncountable, we are finished.
\end{proof}

\section{Proof of Theorem A} \label{Proof of Theorem A}

\subsection{A set of equivalent conditions} \label{ladders}

Analogous to \cite{DV2} we need equivalent versions of the conditions in Theorem A. The obvious proofs are left to the reader.
\begin{lemma} \label{sets} 
Let $(a_{i}),(b_{i})$ be finitely supported sequences of integers, both not identically zero, and put $q_i = a_i-b_i$. The following sets of conditions are equivalent. 
\begin{enumerate} 
\item Let $q_{\mu}$ be the lowest non-zero $q_{i}$ and $q_{\nu}$ the highest non-zero $q_{i}$. 
\begin{enumerate} 
\item $a_{l} = 0$ for $l < \mu$ and $l \geq \nu$.
\item $a_{\mu} = q_{\mu} > 0$.
\item $\sum_{i}q_{i} = 0$
\item $\max(q_l,0) \leq a_l \leq \sum_{i \leq l} q_i$ for all integers $l$.
\end{enumerate} 
\item Let $a_{\mu}$ be the lowest non-zero $a_i$ and $b_{\nu}$ the highest non-zero $b_{i}$. 
\begin{enumerate} 
\item 
The $(a_i),(b_i)$ are non-negative.
\item 
$a_{l} = 0$ for $l \geq \nu$, $b_{l} = 0$ for $l \leq \mu$. 
\item
$\sum_{i}a_{i} = \sum_{i}b_{i}$.
\item
$\sum_{i \le l} b_i \leq \sum_{i < l}a_{i}$ for all integers $l$. 
\end{enumerate} 
\item Put $m = \sum_i a_i$, $n = \sum_i b_i$. 
\begin{enumerate} 
\item The $(a_i),(b_i)$ are non-negative. 
\item $m = n$. 
\item \label{ref-3c-35} $\forall (\alpha, \beta) \in  R_{m,n}: \beta \geq \alpha \Rightarrow (\alpha, \beta)  \in L_{a,b}$. 
\end{enumerate} 
\end{enumerate} 
\end{lemma} 
\begin{lemma} \label{sets2} 
Let $(a_{i}),(b_{i})$ be finitely supported sequences of integers, both not identically zero, and put $q_i = a_i-b_i$. The following sets of conditions are equivalent. 
\begin{enumerate} 
\item Let $q_{\mu}$ be the lowest non-zero $q_{i}$ and $q_{\nu}$ the highest non-zero $q_{i}$. 
\begin{enumerate} 
\item $a_{l} = 0$ for $l < \mu$ and $l \geq \nu$.
\item $a_{\mu} = q_{\mu} > 0$.
\item $\sum_{i}q_{i} = 0$
\item $\max(q_l,0) \leq a_l < \sum_{i \leq l} q_i$ for $\mu < l < \nu$.
\item
If $A$ is cubic it is not true that {\em (}$a_{\mu} \geq 2$ and $\mu = \nu - 1${\em )}.
\end{enumerate} 
\item Let $a_{\mu}$ be the lowest non-zero $a_i$ and $b_{\nu}$ the highest non-zero $b_{i}$. 
\begin{enumerate} 
\item 
The $(a_i),(b_i)$ are non-negative.
\item 
$a_{l} = 0$ for $l \geq \nu$, $b_{l} = 0$ for $l \leq \mu$. 
\item
$\sum_{i}a_{i} = \sum_{i}b_{i}$.
\item
$\sum_{i \le l} b_i < \sum_{i < l}a_{i}$ for $\mu < l < \nu$. 
\item 
If $A$ is cubic it is not true that {\em (}$n \geq 2$ and $\mu = \nu - 1${\em )}.
\end{enumerate} 
\item Put $m = \sum_i a_i$, $n = \sum_i b_i$. 
\begin{enumerate} 
\item The $(a_i),(b_i)$ are non-negative. 
\item $m = n$. 
\item $\forall (\alpha, \beta) \in  R_{m,n}: \beta \geq \alpha - 1 \Rightarrow (\alpha, \beta)  \in L_{a,b}$.  
\item
If $A$ is cubic it is not true that {\em (}$n\geq 2$ and $\forall \alpha, \beta: S(b)_{\beta}- S(a)_{\alpha} = 1${\em )}.
\end{enumerate} 
\end{enumerate} 
\end{lemma} 

\subsection{Proof that the conditions in Theorem A are necessary}
\label{necessary}
This was proved in \cite{Ajitabh} for quadratic $A$, and it is easy to extend it for cubic $A$. As the notations in \cite{Ajitabh} are quite different as in this current paper, we recall the arguments.

\subsubsection{Proof that the conditions in Theorem A(1) are necessary}
\label{Proof that the conditions in Theorem A(1) are necessary}
We will show that the equivalent conditions given in Lemma \ref{sets}(2) are necessary. 

Assume that $M \in \grmod(A)$ is Cohen-Macaulay of GK-dimension two. 
Consider a minimal projective resolution of $M$
\begin{equation} \label{resolutie}
0 \r \bigoplus_i A(-i)^{b_i}\xrightarrow{\phi} \bigoplus_i A(-i)^{a_i}\r M \r 0 
\end{equation}
There is nothing to prove for (2a) and
expressing that $M$ has rank zero gives (2c), so we discuss (2b) and (2d).
The resolution \eqref{resolutie} contains, for all integers $l$, a subcomplex of the form 
\begin{equation} \label{resolutiephil}
\bigoplus_{i \le l} A(-i)^{b_i} \xrightarrow{\phi_l} \bigoplus_{i \leq l} A(-i)^{a_i} 
\end{equation}
Since \eqref{resolutie} is minimal all nonzero entries of a matrix representing $\phi$ have positive degree. Hence the image of $\bigoplus_{i \le l} A(-i)^{b_i}$ under $\phi_{l}$ is contained in $\bigoplus_{i < l} A(-i)^{a_i}$. The fact that $\phi_{l}$ must be injective implies
\begin{equation} \label{ineq}
\sum_{i \le l}{b_i} \leq \sum_{i < l}{a_i}
\end{equation}
from which we obtain (2d).
In particular, if we take $l = \mu$ this shows that 
$b_{i} = 0$ for $i \leq \mu$. 
In order to prove that $a_{i} = 0$ for $i \geq \nu$, add $\sum_{i > l}(a_{i} - b_{i})$ on both sides of \eqref{ineq} and use (2c) to obtain
\[
a_{l} + \sum_{i > l}(a_{i} - b_{i}) \leq \sum_{i}(a_{i} - b_{i}) = 0
\]
thus
\begin{equation} \label{dualineq}
\sum_{l \leq i}{a_i} \leq \sum_{l < i}{b_i}
\end{equation}
Taking $l = \nu$ gives $a_{i} = 0$ for $i \geq \nu$. \\
This completes the proof that the conditions in Theorem A(1) are necessary.

\subsubsection{Proof that the conditions in Theorem A(2) are necessary}
We will show that the equivalent conditions given in Lemma \ref{sets2}(2) are necessary. 

Let $M$ be a critical Cohen-Macaulay module of GK-dimension two. Same reasoning as in \S\ref{Proof that the conditions in Theorem A(1) are necessary} shows Lemma \ref{sets2}(2)(a-c). So we need to show that Lemma \ref{sets2}(2)(d-e) holds.

We will start with the proof of Lemma \ref{sets2}(2)(d), i.e.
\[
\sum_{i \leq l}b_{i} < \sum_{i < l}a_{i} \mbox{ for } \mu < l < \nu
\]
So assume by contradiction that there is some integer $l$ where $\mu < l < \nu$ such that $\sum_{i \le l}{b_i} = \sum_{i < l}{a_i}$. This means that, for the injective map \eqref{resolutiephil}, $\coker \phi_{l}$ has GK-dimension $\leq 2$ and is different from zero.

Note that $\bigoplus_{i < l} A(-i)^{a_i}$ is not zero since $l > \mu$. We have a map $\coker \phi_l \r M$ which we claim to be nonzero. Indeed, if this were the zero map then $\bigoplus_{i < l} A(-i)^{a_i} \r M$ is the zero map, which contradicts the minimality of the resolution \eqref{resolutie}. Hence $\coker \phi_l \r M$ is nonzero. From this we get $\gkdim (\coker \phi_{l}) \leq 2$. Applying $\Hom_A(\coker \phi_l,-)$ to \eqref{resolutie} it is easy to see that actually $\gkdim (\coker \phi_{l}) = 2$.

We will compare the multiplicity $e_{l}$ of $\coker \phi_{l}$ with the multiplicity $e_M$ of $M$. As in the introduction, put $\epsilon = \iota_A e_M$ and $\epsilon_l = \iota_A e_l$. By \eqref{generalHilbertseries} and \eqref{characteristic} we have
\[
\epsilon = \sum_{i}i(b_{i} - a_{i}) \text{ and } \epsilon_{l} = \sum_{i < l}i(b_{i} - a_{i}) + lb_{l} 
\]
We claim that $lb_{l} < \sum_{l \leq i}i(b_{i} - a_{i})$. Indeed, this follows from 
\begin{align*}
\sum_{l \leq i}i(b_{i} - a_{i}) & = 
l\sum_{l \leq i}(b_{i} - a_{i}) + \sum_{l+1 \leq i}(b_{i} - a_{i}) + 
\sum_{l+2 \leq i}(b_{i} - a_{i}) + \ldots \\
& \geq lb_{l} + b_{l+1} + b_{l+2} + \ldots \\
& > lb_{l}
\end{align*}
where the first inequality follows from \eqref{dualineq} and the second one from the assumption that $l < \nu$. Thus we obtain
\begin{align*}
\epsilon_{l} < \sum_{i < l}i(b_{i} - a_{i}) + \sum_{l \leq i}i(b_{i} - a_{i}) = \epsilon_M 
\end{align*}
This means that $\coker \phi_{l}$ has lower multiplicity than $M$. Hence the induced map $\coker \phi_l \r M$ must be zero since $M$ is assumed to be critical. But, as pointed out above, this implies that $\bigoplus_{i < l} A(-i)^{a_i}\r M$ is the zero map, which is impossible. This proves Lemma \ref{sets2}(2d).

What is left to prove is that Lemma \ref{sets2}(2e) holds. If, by contradiction, Lemma \ref{sets2}(2e) is not true then $A$ is cubic and $M$ admits a minimal resolution of the form
\[
0 \r A(-\nu)^n \r A(-(\nu-1))^n \r M \r 0
\]
By shift of grading, we may assume $\nu = 1$. We present the proof for $n = 2$. The arguments are easily extended for all $n \geq 2$. This will complete the proof that the conditions in Theorem A(2) are necessary.
\begin{example} \label{examplestriking}
Assume $A$ is cubic and $M$ is an $A$-module admitting a minimal resolution of the form
\begin{equation} \label{sescubic}
0 \r A(-1)^2 \xrightarrow{
\begin{pmatrix}
l_1 & l_2 \\
l_3 & l_4
\end{pmatrix} \cdot
} A^2 \r M \r 0
\end{equation}
where the entries $l_i = \alpha_i x + \beta_i y \in A_1$ are linear forms ($\alpha_i,\beta_i \in k$). Since
\[
h_M(t) = h_A(t)(2 - 2t) = \frac{2}{(1-t)^2(1+t)} = 2 + 2t + 4t^2 + 4t^3 + 6t^4 + \dots 
\]
we have $\gkdim M = 2$, $e_M = 1$ and $\epsilon_M = 2$. We will show that $M$ is not critical. Let $(x_0,y_0) \in \PP^1$ be a solution of the quadratic equation   
\[
\det 
\begin{pmatrix}
\alpha_1 x_0 + \beta_1 y_0  & \alpha_2 x_0 + \beta_2 y_0 \\
\alpha_3 x_0 + \beta_3 y_0  & \alpha_4 x_0 + \beta_4 y_0
\end{pmatrix} = 0
\]
Thus there is a nonzero $(\lambda, \mu) \in k^2$ for which
\[
\begin{pmatrix}
\alpha_1 x_0 + \beta_1 y_0  & \alpha_2 x_0 + \beta_2 y_0 \\
\alpha_3 x_0 + \beta_3 y_0  & \alpha_4 x_0 + \beta_4 y_0
\end{pmatrix}
\begin{pmatrix}
\la \\
\mu
\end{pmatrix}
 = 0
\]
Consider the linear form $l = y_0 x - x_0 y \in A_1$. Up to scalar multiplication, $l$ is the unique linear form $\alpha x + \beta y$ for which $\alpha x_0 + \beta y_0 = 0$. This means that
\[
\begin{pmatrix}
l_1 & l_2 \\
l_3 & l_4
\end{pmatrix}
\begin{pmatrix}
\la \\ \mu
\end{pmatrix}
= 
\begin{pmatrix}
\gamma \\ \delta 
\end{pmatrix}
l 
\]
for some $\gamma, \delta \in k$. Note that $(\gamma,\delta) \neq (0,0)$ since \eqref{sescubic} is exact. This leads to a commutative diagram
\[
\xymatrix{
0 \ar[r] & A(-1) \ar[d]_{
\begin{pmatrix}
\la \\ \mu
\end{pmatrix}\cdot}
\ar[rr]^{l \cdot} && A \ar[d]^{
\begin{pmatrix}
\gamma \\ \delta
\end{pmatrix}\cdot} \ar[r] & A/lA \ar[r] & 0 \\
0 \ar[r] & A(-1)^2 \ar[rr]^(0.5)
{\begin{pmatrix}
l_1 & l_2 \\
l_3 & l_4
\end{pmatrix}\cdot
}
&& A^2 \ar[r]& M \ar[r] & 0
}
\]
Hence there is a nonzero map $A/lA \r M$. As $A/lA$ has multiplicity $1/2$ and $M$ has multiplicity $1$, this shows that $M$ is not critical, a contradiction.
\end{example}

\subsection{Proof that the conditions in Theorem A(1) are sufficient} \label{Theorem A(1)sufficient}

We fix finitely supported sequences $(a_{i}),(b_{i})$ of non-negative integers such that $\sum_{i}a_{i} = \sum_{i}b_{i} = n$ and we assume the ladder condition holds:
\begin{eqnarray} \label{laddercond}
\forall (\alpha, \beta) \in  R: \beta \geq \alpha \Rightarrow (\alpha, \beta)  \in L_{a,b}
\end{eqnarray}
Thus $S(b)_{\alpha} - S(a)_{\alpha} > 0$ for $1 \leq \alpha \leq n$. Pick nonzero homogeneous elements $h_{\alpha \alpha} \in B_{S(b)_{\alpha} - S(a)_{\alpha}}$. As $A$ is a domain, multiplication by $h_{\alpha \alpha}$ is injective. Let $H_{\alpha}$ be the corresponding cokernels
\[
0 \r A(-S(b)_\alpha) \xrightarrow{h_{\alpha \alpha}\cdot} A(- S(a)_{\alpha}) \r H_\alpha \r 0
\]
for $1 \leq \alpha \leq n$. Then $A$-module $M = H_1 \oplus \dots \oplus H_n$ admits a minimal resolution 
\[
0 \r \oplus_{i}A(-i)^{b_i} \xrightarrow{N \cdot} \oplus_{i}A(-i)^{a_i} \r M \r 0
\]
where
\begin{eqnarray*} 
N =
\begin{pmatrix}
h_{11} & 0 & \dots & 0 \\
0 & h_{22} & \dots & 0 \\
\vdots & & & \vdots \\
0 & 0 & \dots & h_{nn} 
\end{pmatrix}
\end{eqnarray*}
Hence $M$ has projective dimension one, with graded Betti numbers $(a_i)$, $(b_i)$. That $M$ has GK-dimension two is easy to see (see also the proof of Lemma \ref{injective}(3) below). As we have chosen $h_{\alpha \alpha} \in B$, the cyclic modules $H_{\alpha}$ are $g$-torsion free. Hence $M$ is also $g$-torsion free. This completes the proof.

\subsection{Proof that the conditions in Theorem A(2) are sufficient} \label{Proof that the conditions in Theorem A(2) are sufficient}

We fix finitely supported sequences $(a_{i}),(b_{i})$ of non-negative integers for which $\sum_{i}a_{i} = \sum_{i}b_{i} = n$ and we assume that the ladder condition 
\[
\forall (\alpha, \beta) \in  R: \beta \geq \alpha - 1 \Rightarrow (\alpha, \beta)  \in L_{a,b}
\]
is true, together with condition Theorem A(2d):
\[
\text{If $A$ is cubic it is not true that ($n\geq 2$ and $\forall \alpha, \beta: S(b)_{\beta}- S(a)_{\alpha} = 1$)}.
\] 
Put $\epsilon = \sum_i i(b_i - a_i)$. 
We are motivated by the following 
\begin{lemma} \label{injective}
Assume that we have a map $N: \bigoplus_{i}\Oscr_{E}(-i)^{b_{i}} \r \bigoplus_{i}\Oscr_{E}(-i)^{a_{i}}$
such that $X_N(p)$ has 
maximal rank for all but finitely many points $p \in E$. Then 
\begin{enumerate}
\item
$N$ is injective, i.e. we have a short exact sequence in $\coh(E)$
\[
0 \r \bigoplus_{i}\Oscr_{E}(-i)^{b_{i}} \xrightarrow{N} \bigoplus_{i}\Oscr_{E}(-i)^{a_{i}} \r \Nscr \r 0
\]
where $\Nscr \in \coh(E)$ has finite length.
\item
Applying $\Gamma_{\ast}$ to $N$ induces a short exact sequence in $\grmod(B)$
\[
0 \r \bigoplus_{i}B(-i)^{b_{i}} \xrightarrow{\Gamma_{\ast}(N)} \bigoplus_{i}B(-i)^{a_{i}} \r M' \r 0
\]
where $M' \in \grmod(B)$ is pure of GK-dimension one and $\widetilde{M'} = \Nscr$.
\item
Restricting $\Gamma_{\ast}(N)$ to $A$ induces a short exact sequence in $\grmod(A)$
\[
0 \r \bigoplus_{i}A(-i)^{b_{i}} \xrightarrow{\Gamma_{\ast}(N)_A} \bigoplus_{i}A(-i)^{a_{i}} \r M \r 0
\]
where $M \in \grmod(A)$ has GK-dimension two and $\epsilon_M = \sum_{i}i(b_{i} - a_{i})$. Moreover, $M/Mg = M'$ and $M$ is $g$-torsion free.
\end{enumerate}
\end{lemma}
\begin{proof}
(1) If $N$ were not injective then the kernel of $N$ would, as a subsheaf of the vector bundle $\bigoplus_{i}\Oscr_{E}(-i)^{b_{i}}$, have rank $> 0$. Since $\sum_i a_i = \sum_i b_i$ the same is true for the cokernel of $N$. Then $\coker N$ would not be supported on finitely many points in $E$. But this means that $X_N(p)$ has non-maximal rank for infinitely many points $p \in E$, a contradiction.
Thus $N$ is injective. That $\Nscr = \coker N$ has finite length follows from $\sum_i a_i = \sum_i b_i$. 

(2) Apply the functor $\Gamma_{\ast} = \bigoplus_{m \geq 0}H^{0}(E,- \otimes_E \Oscr_E(m))$ to the short exact sequence in (1) and use $\Gamma_{\ast}(\Oscr(l)) = B(l)$ for all integers $l$. As $\Gamma_{\ast}$ is left exact, $\Gamma_{\ast}(N)$ is injective. Since $\underline{\Ext}^{1}_B(k,B) = 0$, $M'$ is not finite dimensional. Hence $M'$ is pure. Application of the exact functor $\widetilde{(-)}$ shows $\widetilde{M'} = \Nscr$ and $\gkdim M' = 1$.
 
(3) The restriction of $\Gamma_{\ast}(N)$ to $A$ defines a map $\bigoplus_{i}A(-i)^{b_{i}} \rightarrow \bigoplus_{i}A(-i)^{a_{i}}$. Write $M = \coker \Gamma_{\ast}(N)_A$. Applying $- \otimes_A B$ and using $\Gamma_{\ast}(N)_A \otimes_A B = \Gamma_{\ast}(N)$, we get $M/Mg = M'$ and $\Tor_1^A(M,B) = 0$. Thus $M$ is $g$-torsion free. As $\gkdim M' = 1$, it is easy to deduce $\gkdim M \leq 2$. Therefore, if the kernel of $\Gamma_{\ast}(N)_A$ is nonzero, it has GK-dimension $\leq 2$. But this is impossible by the pureness of $\bigoplus_{i}A(-i)^{b_{i}}$. 

What remains to prove is $\gkdim M = 2$. By \eqref{characteristic} and Lemma \ref{sets}(2c) we have 
\begin{align*}
2 \lim_{t \r 0} (1-t)^2 h_M(t) & = \sum_{i}i(b_{i} - a_{i}) = \sum_{i}(\nu + 1 - i)(a_{i} - b_{i}) \\
& = \sum_{i \leq \mu}(a_{i} - b_{i}) + 
\sum_{i \leq \mu + 1}(a_{i} - b_{i}) + 
\sum_{i \leq \mu + 2}(a_{i} - b_{i}) + \ldots \\
& \geq \sum_{i}a_{i}
\end{align*}
where we have used Lemma \ref{sets}(2d) to obtain the inequality. Since $\sum_{i}a_{i} > 0$ this proves $\lim_{t \r 0} (1-t)^2 h_M(t) > 0$ i.e. $\gkdim M = 2$, and $e_M = 1/2 \sum_{i}i(b_{i} - a_{i})$. This completes the proof.
\end{proof}
Our proof that the conditions in Theorem A(2) are sufficient follows from the following stronger result.
\begin{theorem} \label{strongthm}
Let $D$ be a multiplicity-free effective divisor of degree $r_A \epsilon - 1$. Then there exists a $g$-torsion free module $M \in \grmod(A)$ of GK-dimension two and projective dimension one which has graded Betti-numbers $(a_i)$, $(b_i)$ i.e. $M$ admits a minimal resolution of the form
\[
0 \r \bigoplus_{i}A(-i)^{b_i} \r \bigoplus_{i}A(-i)^{a_i} \r M \r 0
\]
and for which $\Div(M) = D + (q)$ for some $q \in E$.
\end{theorem}
Indeed, for then we choose a multiplicity-free quantum-irreducible effective divisor $D$ of degree $r_A\epsilon-1$, whose existence is asserted from Theorem \ref{existence}. Lemma \ref{sufficient} implies that the module $M$ in Theorem \ref{strongthm} is critical. 

Thus in order to complete the proof of Theorem A(2) it will be sufficient to prove Theorem \ref{strongthm}. This will be done below.
\begin{proof}[Proof of Theorem \ref{strongthm}] 
Throughout the proof we fix a multiplicity-free effective divisor $D = (q_1) + (q_2) + \dots + (q_{r_A\epsilon-1})$ of degree $r_A\epsilon-1$. As in Lemma \ref{sets2}, let $a_{\mu}$ be the lowest non-zero $a_{i}$ and $b_{\nu}$ be the highest non-zero $b_{i}$. Thus $\mu = S(a)_1$ and $\nu = S(b)_n$. Write $u = \sum_{i < \nu}b_i$ and $v = b_\nu - 1$. 

We break up the proof into six steps.
\begin{step}
Our first step in the proof is to choose a particular $n \times (n-1)$ matrix of the form (only the nonzero entries are indicated)
\begin{eqnarray*} 
H = [H_U \mid H_V] =
\begin{pmatrix}
h_{11} & h_{12} & \dots & h_{1u} & \vline & h_{1,u+1} & \dots & h_{1,n-1} \\
h_{21} & & &  & \vline & & &   \\
& h_{32} & & & \vline & & &  \\
  &        & \ddots & & \vline & & &  \\
 & & & h_{u+1,u} & \vline & & &   \\
  & &   & & \vline & h_{u+2,u+1} &  &   \\
 & &  &  & \vline & & \ddots &  \\
 &  &    &  & \vline &  &  & h_{n,n-1} \\
\end{pmatrix}
\end{eqnarray*} 
whose entries are homogeneous forms $h_{1 \beta} \in B_{S(b)_{\beta} - S(a)_{1}}$, $h_{\beta+1,\beta} \in B_{S(b)_{\beta} - S(a)_{\beta+1}}$ satisfying the following conditions:
\begin{enumerate}
\item[($\ast$)]
The divisors $\Div h_{\alpha \beta}^{S(a)_{\alpha}}$ of the (nonzero) entries $h_{\alpha \beta}$ in $H_U$ are multiplicity-free and have pairwise disjoint support.
\item[($\ast\ast$)]
The support of the divisor $\Div h_{\alpha \beta}^{S(a)_{\alpha}}$ of any entry $h_{\alpha \beta}$ in $H_U$ is disjoint with the support of the divisor $\Div h_{\alpha' \beta'}^{S(a)_{\alpha'}}$ of any entry $h_{\alpha' \beta'}$ in $H_V$. 
\item[($\ast\ast\ast$)]
The divisors $\Div h_{\alpha' \beta'}^{S(a)_{\alpha'}}$ of the entries $h_{\alpha' \beta'}$ in $H_V$ are multiplicity-free. They have pairwise disjoint support unless they appear in the same column $\beta' - u$ of $H_V$. In that case, $\Supp(\Div h_{1,\beta'}^{S(a)_{1}}) \cap \Supp(\Div h_{\beta'+1, \beta'}^{S(a)_{\beta'+1}}) = \{ q_{\beta' - u} \}$. 
\end{enumerate}
Observe that, due to Lemma \ref{multfree}, it is possible to choose such matrices $H_U$, $H_V$ except in the following situation:
\begin{equation} \label{wrong}
\text{$A$ is cubic, $v > 0$ and two linear forms appear in the same column of $H_V$}
\end{equation}
This is because for any two linear forms in $A$ (where $A$ is cubic) their divisors have either disjoint support or the same support (being two distinct points). However, by \eqref{degreematrix} it is easy to see that \eqref{wrong} is same as saying that $n \geq 2$ and all entries of $H$ are linear forms, i.e. $\mu = \nu - 1$. By Lemma \ref{sets2}, this is exactly excluded by condition Theorem A(2d)! In other words, \eqref{wrong} does not occur.
\end{step}
\begin{step}
By construction, the matrix $H$ in Step 1 represents a map 
\[
H: \bigoplus_{i<\nu}\Oscr_E(-i)^{b_i} \oplus \Oscr_E(-\nu)^{b_{\nu}-1} \r \bigoplus_{i}\Oscr_E(-i)^{a_i}
\]
Recall \S\ref{Group law} that in this case the matrix $X_H$ is given by
\begin{eqnarray*} 
X_H =
\begin{pmatrix}
h_{11}^{\sigma^{S(a)_1}} & h_{12}^{\sigma^{S(a)_1}} & \dots & h_{1u}^{\sigma^{S(a)_1}} & \vline & h_{1,u+1}^{\sigma^{S(a)_1}} & \dots & h_{1,n-1}^{\sigma^{S(a)_1}} \\
h_{21}^{\sigma^{S(a)_2}} & & &  & \vline & & &   \\
& h_{32}^{\sigma^{S(a)_3}} & & & \vline & & &  \\
  &        & \ddots & & \vline & & &  \\
 & & & h_{u+1,u}^{\sigma^{S(a)_{u+1}}} & \vline & & &   \\
  & &   & & \vline & h_{u+2,u+1}^{\sigma^{S(a)_{u+2}}} &  &   \\
 & &  &  & \vline & & \ddots &  \\
 &  &    &  & \vline &  &  & h_{n,n-1}^{\sigma^{S(a)_n}} \\
\end{pmatrix}
\end{eqnarray*} 
Therefore, by Step 1 we find
\begin{eqnarray*} 
\rank X_H(p) =
\left\{
\begin{array}{ll}
n-1 & \text{ if } p \in E \setminus \{ q_1,\dots,q_v \} \\
n-2 & \text{ if } p \in \{ q_1,\dots,q_v \}
\end{array}
\right.
\end{eqnarray*}
\end{step}
\begin{step}
Any choice of homogeneous forms $d_\alpha \in B_{S(b)_n - S(a)_\alpha}$, $\alpha = 1, \dots, n$ determines a matrix 
\begin{eqnarray}
[H \mid d] = 
\begin{pmatrix}
h_{11} & h_{12} & \dots  & h_{1,n-1} & \vline & d_1 \\
h_{21} &        &        &           & \vline & d_2 \\
       & h_{32} &        &           & \vline & d_3 \\
       &        & \ddots &           & \vline & \vdots \\ 
       &        & \      & h_{n,n-1} & \vline & d_{n} 
\end{pmatrix}
\end{eqnarray}
which represents a map $[H \mid d]: \bigoplus_{i}\Oscr_C(-i)^{b_i} \r \bigoplus_{i}\Oscr_C(-i)^{a_i}$. We then consider the $k$-linear map
\begin{align*}
\theta: \bigoplus_{\alpha = 1}^{n}B_{S(b)_n - S(a)_\alpha} \r H^{0}(E,\Lscr^{p(\sigma)}): d \mapsto 
\det X_{[H \mid d]}
\end{align*}
where $p(t) = \sum_i(a_i - b_i)t^i/(1-t) \in \ZZ[t,t^{-1}]$. Furthermore, by Step 2 the image of $\theta$ is contained in the $k$-linear subspace
\[
W = \{ s \in H^{0}(E,\Lscr^{p(\sigma)}) \mid s(q_i) = 0 \text{ for } i = 1,\dots,v \} \subset H^{0}(E,\Lscr^{p(\sigma)})
\] 
We claim that $\im \theta = W$. This will follow from the Steps 4 and 5 below.
\end{step}
\begin{step}
If $\dim_k \ker \theta = \sum_{\alpha = 1}^{n-1} \dim_k B_{S(b)_n - S(b)_\alpha}$ then $\im \theta = W$. Indeed, (a generalized version of) Lemma \ref{dimdiv} shows that $\codim W = v$. Thus
\begin{align*}
\dim_k W = \dim_k H^{0}(E,\Lscr^{p(\sigma)}) - v & = \sum_{\alpha=1}^n \dim_k B_{S(b)_\alpha - S(a)_\alpha} - v \\
& = \sum_{\alpha=1}^n r_A(S(b)_\alpha -S(a)_\alpha) - v 
\end{align*}
while on the other hand
\begin{align*}
\dim_k \im \theta & = \sum_{\alpha = 1}^n \dim_k B_{S(b)_n - S(a)_\alpha} - \dim_k \ker \theta \\
& = \sum_{\alpha = 1}^n \dim_k B_{S(b)_n - S(a)_\alpha} - \sum_{\alpha = 1}^{n-1} \dim_k B_{S(b)_n - S(b)_\alpha} \\
& = \sum_{\alpha= 1}^{n}r_A(S(b)_n - S(a)_\alpha) - \sum_{\alpha = 1}^{n}r_A(S(b)_n - S(b)_\alpha) - (b_\tau - 1) \\
& = \sum_{\alpha=1}^n r_A(S(b)_\alpha -S(a)_\alpha) - v
\end{align*}
\end{step}
\begin{step}
$\dim_k \ker \theta = \sum_{\alpha = 1}^{n-1} \dim_k B_{S(b)_n - S(b)_\alpha}$. We prove this as follows. For any choice of homogeneous elements $c_\alpha \in B_{S(b)_n - S(b)_\alpha}$ for $\alpha = 1,\dots,n-1$, putting
\begin{eqnarray} \label{formkernel}
\begin{pmatrix}
d_1 \\
d_2 \\
d_3 \\
\vdots \\
d_n
\end{pmatrix}
=
\begin{pmatrix}
h_{11} & h_{12} & \dots  & h_{1,n-1} \\
h_{21} &        &        &           \\
       & h_{32} &        &           \\
       &        & \ddots &           \\ 
       &        & \      & h_{n,n-1}  
\end{pmatrix}
\begin{pmatrix}
c_1 \\
c_2 \\
\vdots \\
c_{n-1} 
\end{pmatrix}
\end{eqnarray}
yields an element $d = (d_1, \dots, d_n) \in \bigoplus_{\alpha = 1}^{n}B_{S(b)_n - S(a)_\alpha}$ in the kernel of $\theta$. 
%
%
Thus we have a $k$-linear map
\begin{align*}
\tilde{\theta}: \bigoplus_{\alpha = 1}^{n-1}B_{S(b)_n - S(b)_\alpha} \r \ker \theta: (c_1,\dots,c_{n-1})^t \mapsto H \cdot (c_1,\dots,c_{n-1})^t
\end{align*}
which is injective by the fact that the entries of $H$ are nonzero (Step 1) and $B$ is a domain. Hence in order to prove Step 5 it suffices to show $\tilde{\theta}$ is surjective.

Pick $d = (d_1, \dots, d_n) \in \ker \theta$. By Step 2 we may solve \eqref{formkernel} locally at $p \in E\setminus \{q_1, \dots,q_{v}\}$, i.e. we may find a solution $\la(p) = (\la_1(p), \dots, \la_{n-1}(p))$, where 
\[
\la_\alpha(p) \in \left( \Lscr_{S(b)_n - S(b)_\alpha}^{\sigma^{S(b)_{\alpha}}} \right)_p / m_p \left( \Lscr_{S(b)_n - S(b)_\alpha}^{\sigma^{S(b)_{\alpha}}} \right)_p 
\]
such that
\begin{multline} \label{equivkernel}
\begin{pmatrix}
d_1^{\sigma^{S(a)_1}}(p) \\
d_2^{\sigma^{S(a)_2}}(p) \\
d_3^{\sigma^{S(a)_3}}(p) \\
\vdots \\
d_n^{\sigma^{S(a)_n}}(p)
\end{pmatrix}
=
\begin{pmatrix}
h_{11}^{\sigma^{S(a)_1}}(p) & h_{12}^{\sigma^{S(a)_1}}(p) & \dots  & h_{1,n-1}^{\sigma^{S(a)_1}}(p) \\
h_{21}^{\sigma^{S(a)_2}}(p) &        &        &           \\
       & h_{32}^{\sigma^{S(a)_3}}(p) &        &           \\
       &        & \ddots &           \\ 
       &        & \      & h_{n,n-1}^{\sigma^{S(a)_n}}(p)  
\end{pmatrix}
\otimes
\begin{pmatrix}
\la_1(p) \\
\la_2(p) \\
\vdots \\
\la_{n-1}(p) 
\end{pmatrix}
\end{multline}
%
%
To show that we can solve \eqref{formkernel} globally, we proceed as follows.
\begin{itemize}
\item
For $\beta = 1,\dots,u$ the $\beta + 1$-th equation in \eqref{equivkernel} becomes
\begin{equation} \label{eqnh}
d_{\beta+1}(p^{\sigma^{S(a)_{\beta+1}}}) 
= h_{\beta+1,\beta}(p^{\sigma^{S(a)_{\beta+1}}}) \otimes \la_{\beta}(p) \text{ for all } p \in E\setminus \{q_1, \dots,q_{v}\}
\end{equation}
By Step 1, $q_1,\dots,q_v \not \in \Div h_{\beta+1,\beta}^{\sigma^{S(a)_{\beta+1}}}$. 
Hence we deduce from \eqref{eqnh}
\[
d_{\beta+1} = 0 \text{ or } \Div d_{\beta+1} = \Div h_{\beta+1,\beta} + D'
\]
for some effective divisor $D'$. By Lemma \ref{division} this means that $d_{\beta+1} = h_{\beta+1,\beta}c_\beta$ for some $c_\beta \in B_{S(b)_n - S(b)_\beta}$. 
\item
For $\beta = u+1,\dots,n-1$ the $\beta + 1$-th equation in \eqref{equivkernel} becomes
\begin{equation} \label{eqnh2}
d_{\beta+1}(p^{\sigma^{S(a)_{\beta+1}}})
= h_{\beta+1,\beta}(p^{\sigma^{S(a)_{\beta+1}}})\otimes \la_{\beta}(p) \text{ for all } p \in E\setminus \{q_1, \dots,q_{v}\}
\end{equation}
As $S(b)_{\beta} = S(b)_n$ we have $\deg d_{\beta+1} = \deg h_{\beta+1,\beta}$ (if $d_{\beta+1} \neq 0$). By Step 1 there is only one $i = 1,\dots,v$ for which $q_i^{\sigma^{S(a)_{\beta+1}}} \in \Div h_{\beta+1,\beta}$. As $h_{\beta+1,\beta}^{S(a)_{\beta+1}}$ is multiplicity-free Lemma \ref{division'} yields $d_{\beta+1} = h_{\beta+1,\beta}c_\beta$ for some $c_\beta \in B_{S(b)_n - S(b)_\beta} = k$. 
\item
Finally, by the previous two items the first equation in \eqref{equivkernel} becomes
\begin{align*}
d_1(p^{\sigma^{S(a)_1}})
& = h_{11}(p^{\sigma^{S(a)_1}})\otimes \la_1(p) +  
\dots + h_{1,n-1}(p^{\sigma^{S(a)_1}})\otimes\la_{n-1}(p) \\
& = (h_{11}c_1 + h_{12}c_2 + \dots + h_{1,n-1}c_{n-1})(p^{\sigma^{S(a)_1}})
\end{align*}
for $p \in E\setminus \{ q_1,\dots,q_v \}$.
Hence $(d_1 - \sum_{\beta = 1}^{n-1}h_{1 \beta}c_{\beta})(p) = 0$ for all but finitely many $p \in E$. This clearly implies $d_1 = \sum_{\beta = 1}^{n-1}h_{1 \beta}c_{\beta}$. 
\end{itemize}
It follows that $d$ is of the form \eqref{formkernel}. We have shown that $\tilde{\theta}$ is surjective. This ends the proof of Step 5.
\end{step}
\begin{step}
As $\dim_k H^0(E,\Lscr^{p(\sigma)}) = r_A \epsilon$, we may pick a global section $s \in H^0(E,\Lscr^{\sigma})$ for which $\Div(s) = D + (q)$ for some $q \in E$. Clearly $s \in W$. By Steps 4 and 5, we have $\im \theta = W$. Thus we may find homogeneous forms $d_\alpha \in B_{S(b)_n - S(a)_\alpha}$, $\alpha = 1, \dots, n$ for which $\det X_{[H \mid d]} = s$. By Lemma \ref{injective}, there is an short exact sequence in $\grmod(A)$
\[
0 \r \bigoplus_{i}A(-i)^{b_{i}} \xrightarrow{\Gamma_{\ast}([H \mid d])_A} \bigoplus_{i}A(-i)^{a_{i}} \r M \r 0
\]
where $M$ is $g$-torsion free of GK-dimension two. By construction, $\Div(M) = \Div(s) = D + (p)$. This completes the proof of Theorem \ref{strongthm}.
\qed
\end{step}
\def\qed{}
\end{proof}

\section{Proof of Theorem B and other properties of Hilbert series}
\label{Proof of Theorem B and other properties of Hilbert series}

\begin{proof}[Proof of Theorem B]
First, let $M$ be a normalized Cohen-Macaulay $A$-module of GK-dimension two and multiplicity $e$. Writing the Hilbert series $h_M(t)$ of $M$ in the form \eqref{generalHilbertseries} we see that there is a Laurent polynomial $s(t)$ for which 
\begin{equation} \label{maybewritten}
h_M(t) = h_A(t)(\epsilon(1-t) - s(t)(1-t)^2)
\end{equation}
where $\epsilon = \iota_A e$. Since $M$ is normalized we have $M_{< 0} = 0$, thus $s(t) \in \ZZ[t]$. 

Let $(a_i)$, $(b_i)$ denote the graded Betti numbers of $M$ and consider the characteristic polynomial $q_M(t) = \sum_i(a_i - b_i)t^i$. Then $q_M(t)/(1-t) = \sum_{l} p_l t^l$ where $p_l = \sum_{i \leq l}q_i$. By \S\ref{Theorem A(1)sufficient}, the conditions of Lemma \ref{sets}(1)(a-d) hold. Note that, as $M$ is normal, $a_0$ is the lowest non-zero $a_i$ i.e. $\mu = 0$. In particular, 
\begin{eqnarray*}
p_l 
\left\{
\begin{array}{ll}
> 0 & \text{ for } l = 0 \\
\geq 0 & \text{ for } 0 < l < \nu \\
= 0 & \text{ else }
\end{array}
\right.
\end{eqnarray*}
Combining \eqref{characteristic} and \eqref{maybewritten} we have
\begin{equation} \label{multiply}
s(t)(1-t) = \epsilon - \sum_{l}p_{l}t^{l}
\end{equation}
Multiplying \eqref{multiply} by $1/(1-t) = 1 + t + t^{2} + \dots$ shows that $s(t)$ is of the form
\begin{equation*} 
\epsilon > s_{0} \geq s_{1} \geq \dots \geq 0
\end{equation*}
If $M$ is in addition critical, Lemma \ref{sets2}(1d) implies $p_l > 0$ for $0 \leq l < \nu$. By same reasoning as above we find that $s(t)$ is of the form
\[
\epsilon > s_{0} > s_{1} > \dots \geq 0
\]
In case $A$ is cubic, Lemma \ref{sets2}(1e) requires in addition that $q(t)$ is not of the form $n(1-t)$ for $n = \sum_i a_i = \sum_i b_i \geq 2$. This is the same as saying that in case $\epsilon \geq 2$ $q(t) \neq \epsilon(1-t)$ . In other words, in case $\epsilon \geq 2$ we have $s(t) \neq 0$. 

The converse of Theorem B follows by reversing the arguments.
\end{proof}
%
%
%
%
\begin{remark}
From Theorem B we may deduce the following combinatorical result. For positive integers $m,n$ let $p(D,n,<m)$ denote the number of partitions of $n$ with distinct parts in which every part is strictly smaller than $m$. Needless to say that $p(D,n,<m) = 0$ for $n > m(m-1)/2$. Corollary \ref{count} now yields
\[
\sum_{n \geq 0}p(D,n,<m) = 2^{m-1}
\]
for all positive integers $m$.
\end{remark}
We also mention
\begin{corollary} \label{countresolutions}
Let $\epsilon > 0$ be an integer. The number of finitely supported sequences $(a_i)$, $(b_i)$ which occur as the graded Betti numbers of a (resp. critical) normalized Cohen-Macaulay $A$-module $M$ of GK-dimension two having Hilbert series
\begin{equation*} 
h_M(t) = h_A(t)(\epsilon(1-t) - s(t)(1-t)^2)
\end{equation*}
is equal to 
\[
[1 + \min (\epsilon - s_{0}, s_{0} - s_{1})] \cdot \prod_{1 < l} [1 + \min (s_{l-2} - s_{l-1}, s_{l-1} - s_{l})]
\]
resp.
\begin{equation} \label{number}
\min (\epsilon - s_{0}, s_{0} - s_{1}) \cdot \prod_{1 < l} \min(s_{l-2} - s_{l-1}, s_{l-1} - s_{l})
\end{equation}
This number \eqref{number} is bigger than one if and only if there are two consecutive downward jumps of length $\geq 2$ in the coefficients of $\epsilon t^{-1} + s(t)$. 
\end{corollary} 
\begin{proof}
The number of solutions to the conditions Lemma \ref{sets}(1)(a-d) is
\[ 
\prod_{\mu < l < \nu} \left(\biggl(\sum_{i \leq l}q_i \biggr) - \max(q_l,0) + 1 \right) = 
\prod_{l>\mu}\min\biggl( 1 + \sum_{i< l}q_i, 1 + 
\sum_{i\le l}q_i\biggr) 
\] 
Since we restict to normalized modules we have $\mu = 0$.
Noting that $q_{0} = \epsilon - s_{0}$ and $\sum_{i \leq l}q_{i} = s_{l-1} - s_{l}$ for $l > 0$ yields that the number of solutions is equal to
\[
\min (1+ \epsilon - s_{0}, 1 + s_{0} - s_{1}) \cdot \prod_{1 < l} \min( 1 + s_{l-2} - s_{l-1}, 1 + s_{l-1} - s_{l})
\]
Same reasoning in the critical case.
\end{proof}

\appendix 
\section{Hilbert series up to $\epsilon = 4$} \label{A} 
Let $A$ be a generic three-dimensional Artin-Schelter regular algebra, either quadratic or cubic \S\ref{Elliptic algebras}. Let $M$ be a normalized critical Cohen-Macaulay graded right $A$-module of GK-dimension two. According to Theorem B the Hilbert series of $M$ has the form
\[
h_M(t) = h_A(t)(\epsilon(1-t) - s_M(t)(1-t)^2) 
\]
where $\epsilon > 0$ is an integer and $s_M(t)\in \ZZ[t]$ is a polynomial of the form
\begin{equation*} 
\epsilon > s_{0} > s_{1} > \dots \geq 0
\text{ and if $A$ is cubic and $\epsilon > 1$ then $s(t) \neq 0$ }
\end{equation*}
The multiplicity of $M$ is given by $e_M = \epsilon /2$. For the cases $\epsilon \leq 4$ we list the possible Hilbert series for $M$, the corresponding $s(t)$ and the possible minimal resolutions of $M$. Recall that 
\begin{eqnarray*}
r_A = 
\left\{
\begin{array}{ll}
3 & \text{ if $A$ is quadratic} \\
2 & \text{ if $A$ is cubic}
\end{array}
\right.
\end{eqnarray*}

\def\mystrut{\vrule width 0em height 2em}
\[ 
\begin{array}{|c|l|} 
\hline 
\epsilon = 1 
& h_{M}(t)  = 
\left\{
\begin{array}{ll}
1 + 2t + 3t^{2} + 4t^{3} + 5t^{4} + 6t^{5} + \ldots & \text{ if } r_A = 3 \\
1 + t + 2t^2 + 2t^3 + 3t^{4} + 3t^{5} + \ldots & \text{ if } r_A = 2 
\end{array}
\right. \mystrut \\
& s_{M}(t)  = 0  \\
& 0 \r A(-1) \r A \r M \r 0 
\mystrut \\ 
\hline 
%
%
%
\epsilon = 2 
& h_{M}(t)  = 
\left\{
\begin{array}{ll}
2 + 4t + 6t^{2} + 8t^{3} + 10t^{4} + 12t^{5} + \ldots & \text{ if } r_A = 3 \\
\emptyset & \text{ if } r_A = 2 
\end{array}
\right. \mystrut \\
& s_{M}(t)  = 0 \\ 
& 0 \r A(-1)^{2} \r A^{2} \r M \r 0 \\ 
\cline{2-2} 
& h_{M}(t)  = 
\left\{
\begin{array}{ll}
1 + 3t + 5t^{2} + 7t^{3} + 9t^{4} + 11t^{5} + \ldots & \text{ if } r_A = 3 \\
1 + 2t + 3t^2 + 4t^3 + 5t^{4} + 6t^{5} + \ldots & \text{ if } r_A = 2 
\end{array}
\right. \mystrut \\
\smash{\parbox[c]{0.5cm}{\begin{center} 
\unitlength 1mm
\begin{picture}(5.00,10.00)(0,0)
\linethickness{0.15mm}
\put(0.00,5.00){\line(1,0){5.00}}
\put(0.00,5.00){\line(0,1){5.00}}
\put(5.00,5.00){\line(0,1){5.00}}
\put(0.00,10.00){\line(1,0){5.00}}
\end{picture}
\end{center}}} 
& s_{M}(t)  = 1 \\ 
& 0 \r A(-2) \r A \r M \r 0 \\ 
\hline
\epsilon = 3
& h_{M}(t)  = 
\left\{
\begin{array}{ll}
3 + 6t + 9t^{2} + 12t^{3} + 15t^{4} + 18t^{5} + \ldots & \text{ if } r_A = 3  \\
\emptyset & \text{ if } r_A = 2 
\end{array}
\right. \mystrut \\
& s_{M}(t)  = 0 \\ 
& 0 \r A(-1)^{3} \r A^{3} \r M \r 0 \\
\cline{2-2} 
& h_{M}(t)  = 
\left\{
\begin{array}{ll}
2 + 5t + 8t^{2} + 11t^{3} + 14t^{4} + 17t^{5} + \ldots & \text{ if } r_A = 3 \\
2 + 3t + 5t^2 + 6t^3 + 8t^{4} + 9t^{5} + \ldots & \text{ if } r_A = 2 
\end{array}
\right. \mystrut \\
\smash{\parbox[c]{0.5cm}{\begin{center}
\unitlength 1mm
\begin{picture}(5.00,10.00)(0,0)
\linethickness{0.15mm}
\put(0.00,5.00){\line(1,0){5.00}}
\put(0.00,5.00){\line(0,1){5.00}}
\put(5.00,5.00){\line(0,1){5.00}}
\put(0.00,10.00){\line(1,0){5.00}}
\end{picture}
\end{center}}} 
& s_{M}(t)  = 1 \\ 
& 0 \r A(-1) \oplus A(-2) \r A^{2} \r M \r 0 \\ 
\cline{2-2} 
& h_{M}(t)  = 
\left\{
\begin{array}{ll}
1 + 4t + 7t^{2} + 10t^{3} + 13t^{4} + 16t^{5} + \ldots & \text{ if } r_A = 3 \\
1 + 3t + 4t^2 + 6t^3 + 7t^{4} + 9t^{5} + \ldots & \text{ if } r_A = 2 
\end{array}
\right. \mystrut \\
\smash{\parbox[c]{0.5cm}{\begin{center}
\unitlength 1mm
\begin{picture}(5.00,15.00)(0,0)
\linethickness{0.15mm}
\put(0.00,5.00){\line(1,0){5.00}}
\put(0.00,5.00){\line(0,1){5.00}}
\put(5.00,5.00){\line(0,1){5.00}}
\put(0.00,10.00){\line(1,0){5.00}}
\linethickness{0.15mm}
\put(0.00,10.00){\line(1,0){5.00}}
\put(0.00,10.00){\line(0,1){5.00}}
\put(5.00,10.00){\line(0,1){5.00}}
\put(0.00,15.00){\line(1,0){5.00}}
\end{picture}
\end{center}}} 
& s_{M}(t)  = 2 \\ 
& 0 \r A(-2)^{2} \r A \oplus A(-1) \r M \r 0 \\ 
\cline{2-2}
& h_M(t) = 
\left\{
\begin{array}{ll}
1 + 3t + 6t^{2} + 9t^{3} + 12t^{4} + 15t^{5} + \ldots & \text{ if } r_A = 3 \\
1 + 2t + 4t^2 + 5t^3 + 7t^{4} + 8t^{5} + \ldots & \text{ if } r_A = 2 
\end{array}
\right. \mystrut \\ 
\smash{\parbox[c]{1cm}{\begin{center}
\unitlength 1mm
\begin{picture}(10.00,15.00)(0,0)
\linethickness{0.15mm}
\put(0.00,5.00){\line(1,0){5.00}}
\put(0.00,5.00){\line(0,1){5.00}}
\put(5.00,5.00){\line(0,1){5.00}}
\put(0.00,10.00){\line(1,0){5.00}}
\linethickness{0.15mm}
\put(0.00,10.00){\line(1,0){5.00}}
\put(0.00,10.00){\line(0,1){5.00}}
\put(5.00,10.00){\line(0,1){5.00}}
\put(0.00,15.00){\line(1,0){5.00}}
\linethickness{0.15mm}
\put(5.00,5.00){\line(1,0){5.00}}
\put(5.00,5.00){\line(0,1){5.00}}
\put(10.00,5.00){\line(0,1){5.00}}
\put(5.00,10.00){\line(1,0){5.00}}
\end{picture}
\end{center}}} 
& s_{M}(t)  = 2 + t \\ 
& 0 \r A(-3) \r A \r M \r 0 \\ 
\hline 
\end{array}
\]
\[ 
\begin{array}{|c|l|} 
\hline
\epsilon = 4
& h_{M}(t)  = 
\left\{
\begin{array}{ll}
4 + 8t + 12t^{2} + 16t^{3} + 20t^{4} + 24t^{5} + \ldots & \text{ if } r_A = 3  \\
\emptyset & \text{ if } r_A = 2 
\end{array}
\right. \mystrut \\
& s_{M}(t)  = 0 \\ 
& 0 \r A(-1)^{4} \r A^{4} \r M \r 0 \\
\cline{2-2} 
& h_M(t) = 
\left\{
\begin{array}{ll}
3 + 7t + 11t^{2} + 15t^{3} + 19t^{4} + 23t^{5} + \ldots & \text{ if } r_A = 3 \\
3 + 4t + 7t^2 + 8t^3 + 11t^{4} + 12t^{5} + \ldots & \text{ if } r_A = 2 
\end{array}
\right. \mystrut \\
\smash{\parbox[c]{0.5cm}{\begin{center}
\unitlength 1mm
\begin{picture}(5.00,10.00)(0,0)
\linethickness{0.15mm}
\put(0.00,5.00){\line(1,0){5.00}}
\put(0.00,5.00){\line(0,1){5.00}}
\put(5.00,5.00){\line(0,1){5.00}}
\put(0.00,10.00){\line(1,0){5.00}}
\end{picture}
\end{center}}} 
& s_{M}(t)  = 1 \\ 
& 0 \r A(-1)^{2} \oplus A(-2) \r A^{3} \r M \r 0 \\ 
\cline{2-2} 
& h_M(t) = 
\left\{
\begin{array}{ll}
2 + 6t + 10t^{2} + 14t^{3} + 18t^{4} + 22t^{5} + \ldots & \text{ if } r_A = 3 \\
2 + 4t + 6t^2 + 8t^3 + 10t^{4} + 12t^{5} + \ldots & \text{ if } r_A = 2 
\end{array}
\right. \mystrut \\
\smash{\parbox[c]{0.5cm}{\begin{center}
\unitlength 1mm
\begin{picture}(5.00,15.00)(0,0)
\linethickness{0.15mm}
\put(0.00,5.00){\line(1,0){5.00}}
\put(0.00,5.00){\line(0,1){5.00}}
\put(5.00,5.00){\line(0,1){5.00}}
\put(0.00,10.00){\line(1,0){5.00}}
\linethickness{0.15mm}
\put(0.00,10.00){\line(1,0){5.00}}
\put(0.00,10.00){\line(0,1){5.00}}
\put(5.00,10.00){\line(0,1){5.00}}
\put(0.00,15.00){\line(1,0){5.00}}
\end{picture}
\end{center}}} 
& s_{M}(t)  = 2 \\ 
& 0 \r A(-2)^{2} \r A^{2} \r M \r 0 \\ 
& 0 \r A(-1) \oplus A(-2)^{2} \r A^{2} \oplus A(-1) \r M \r 0 \\
\cline{2-2} 
& h_M(t) = 
\left\{
\begin{array}{ll}
2 + 5t + 9t^{2} + 13t^{3} + 17t^{4} + 21t^{5} + \ldots & \text{ if } r_A = 3 \\
2 + 3t + 6t^2 + 7t^3 + 10t^{4} + 11t^{5} + \ldots & \text{ if } r_A = 2 
\end{array}
\right. \mystrut \\ 
\smash{\parbox[c]{1cm}{\begin{center}
\unitlength 1mm
\begin{picture}(10.00,15.00)(0,0)
\linethickness{0.15mm}
\put(0.00,5.00){\line(1,0){5.00}}
\put(0.00,5.00){\line(0,1){5.00}}
\put(5.00,5.00){\line(0,1){5.00}}
\put(0.00,10.00){\line(1,0){5.00}}
\linethickness{0.15mm}
\put(0.00,10.00){\line(1,0){5.00}}
\put(0.00,10.00){\line(0,1){5.00}}
\put(5.00,10.00){\line(0,1){5.00}}
\put(0.00,15.00){\line(1,0){5.00}}
\linethickness{0.15mm}
\put(5.00,5.00){\line(1,0){5.00}}
\put(5.00,5.00){\line(0,1){5.00}}
\put(10.00,5.00){\line(0,1){5.00}}
\put(5.00,10.00){\line(1,0){5.00}}
\end{picture}
\end{center}}} 
& s_{M}(t)  = 2 + t \\ 
& 0 \r A(-1) \oplus A(-3) \r A^{2} \r M \r 0 \\ 
\cline{2-2}
& h_M(t) = 
\left\{
\begin{array}{ll}
1 + 5t + 9t^{2} + 13t^{3} + 17t^{4} + 21t^{5} + \ldots & \text{ if } r_A = 3 \\
1 + 4t + 5t^2 + 8t^3 + 9t^{4} + 12t^{5} + \ldots & \text{ if } r_A = 2 
\end{array}
\right. \mystrut \\
\smash{\parbox[c]{0.5cm}{\begin{center}
\unitlength 1mm
\begin{picture}(5.00,20.00)(0,0)
\linethickness{0.15mm}
\put(0.00,5.00){\line(1,0){5.00}}
\put(0.00,5.00){\line(0,1){5.00}}
\put(5.00,5.00){\line(0,1){5.00}}
\put(0.00,10.00){\line(1,0){5.00}}
\linethickness{0.15mm}
\put(0.00,10.00){\line(1,0){5.00}}
\put(0.00,10.00){\line(0,1){5.00}}
\put(5.00,10.00){\line(0,1){5.00}}
\put(0.00,15.00){\line(1,0){5.00}}
\linethickness{0.15mm}
\put(0.00,15.00){\line(1,0){5.00}}
\put(0.00,15.00){\line(0,1){5.00}}
\put(5.00,15.00){\line(0,1){5.00}}
\put(0.00,20.00){\line(1,0){5.00}}
\end{picture}
\end{center}}}
& s_{M}(t)  = 3 \\ 
& 0 \r A(-2)^{3} \r A \oplus A(-1)^{2} \r M \r 0 \\
\cline{2-2} 
& h_M(t) = 
\left\{
\begin{array}{ll}
1 + 4t + 8t^{2} + 12t^{3} + 16t^{4} + 20t^{5} + \ldots & \text{ if } r_A = 3 \\
1 + 3t + 5t^2 + 7t^3 + 9t^{4} + 11t^{5} + \ldots & \text{ if } r_A = 2 
\end{array}
\right. \mystrut \\
\smash{\parbox[c]{1cm}{\begin{center}
\unitlength 1mm
\begin{picture}(10.00,20.00)(0,0)
\linethickness{0.15mm}
\put(0.00,5.00){\line(1,0){5.00}}
\put(0.00,5.00){\line(0,1){5.00}}
\put(5.00,5.00){\line(0,1){5.00}}
\put(0.00,10.00){\line(1,0){5.00}}
\linethickness{0.15mm}
\put(0.00,10.00){\line(1,0){5.00}}
\put(0.00,10.00){\line(0,1){5.00}}
\put(5.00,10.00){\line(0,1){5.00}}
\put(0.00,15.00){\line(1,0){5.00}}
\linethickness{0.15mm}
\put(0.00,15.00){\line(1,0){5.00}}
\put(0.00,15.00){\line(0,1){5.00}}
\put(5.00,15.00){\line(0,1){5.00}}
\put(0.00,20.00){\line(1,0){5.00}}
\linethickness{0.15mm}
\put(5.00,5.00){\line(1,0){5.00}}
\put(5.00,5.00){\line(0,1){5.00}}
\put(10.00,5.00){\line(0,1){5.00}}
\put(5.00,10.00){\line(1,0){5.00}}
\end{picture}
\end{center}}} 
& s_{M}(t)  = 3 + t \\ 
& 0 \r A(-2) \oplus A(-3) \r A \oplus A(-1) \r M \r 0 \\ 
\cline{2-2} 
& h_M(t) = 
\left\{
\begin{array}{ll}
1 + 3t + 7t^{2} + 11t^{3} + 15t^{4} + 19t^{5} + \ldots & \text{ if } r_A = 3 \\
1 + 2t + 5t^2 + 6t^3 + 9t^{4} + 10t^{5} + \ldots & \text{ if } r_A = 2 
\end{array}
\right. \mystrut \\
\smash{\parbox[c]{1cm}{\begin{center}
\unitlength 1mm
\begin{picture}(10.00,20.00)(0,0)
\linethickness{0.15mm}
\put(0.00,5.00){\line(1,0){5.00}}
\put(0.00,5.00){\line(0,1){5.00}}
\put(5.00,5.00){\line(0,1){5.00}}
\put(0.00,10.00){\line(1,0){5.00}}
\linethickness{0.15mm}
\put(0.00,10.00){\line(1,0){5.00}}
\put(0.00,10.00){\line(0,1){5.00}}
\put(5.00,10.00){\line(0,1){5.00}}
\put(0.00,15.00){\line(1,0){5.00}}
\linethickness{0.15mm}
\put(0.00,15.00){\line(1,0){5.00}}
\put(0.00,15.00){\line(0,1){5.00}}
\put(5.00,15.00){\line(0,1){5.00}}
\put(0.00,20.00){\line(1,0){5.00}}
\linethickness{0.15mm}
\put(5.00,5.00){\line(1,0){5.00}}
\put(5.00,5.00){\line(0,1){5.00}}
\put(10.00,5.00){\line(0,1){5.00}}
\put(5.00,10.00){\line(1,0){5.00}}
\linethickness{0.15mm}
\put(5.00,10.00){\line(1,0){5.00}}
\put(5.00,10.00){\line(0,1){5.00}}
\put(10.00,10.00){\line(0,1){5.00}}
\put(5.00,15.00){\line(1,0){5.00}}
\end{picture}
\end{center}}} 
& s_{M}(t)  = 3 + 2t \\ 
& 0 \r A(-3)^{2} \r A \oplus A(-2) \r M \r 0 \\ 
\cline{2-2}
& h_M(t) = 
\left\{
\begin{array}{ll}
1 + 3t + 6t^{2} + 10t^{3} + 14t^{4} + 18t^{5} + \ldots & \text{ if } r_A = 3 \\
1 + 2t + 4t^2 + 6t^3 + 8t^{4} + 10t^{5} + \ldots & \text{ if } r_A = 2 
\end{array}
\right. \mystrut \\
\smash{\parbox[c]{1.5cm}{\begin{center}
\unitlength 1mm
\begin{picture}(15.00,20.00)(0,0)
\linethickness{0.15mm}
\put(0.00,5.00){\line(1,0){5.00}}
\put(0.00,5.00){\line(0,1){5.00}}
\put(5.00,5.00){\line(0,1){5.00}}
\put(0.00,10.00){\line(1,0){5.00}}
\linethickness{0.15mm}
\put(0.00,10.00){\line(1,0){5.00}}
\put(0.00,10.00){\line(0,1){5.00}}
\put(5.00,10.00){\line(0,1){5.00}}
\put(0.00,15.00){\line(1,0){5.00}}
\linethickness{0.15mm}
\put(0.00,15.00){\line(1,0){5.00}}
\put(0.00,15.00){\line(0,1){5.00}}
\put(5.00,15.00){\line(0,1){5.00}}
\put(0.00,20.00){\line(1,0){5.00}}
\linethickness{0.15mm}
\put(5.00,5.00){\line(1,0){5.00}}
\put(5.00,5.00){\line(0,1){5.00}}
\put(10.00,5.00){\line(0,1){5.00}}
\put(5.00,10.00){\line(1,0){5.00}}
\linethickness{0.15mm}
\put(5.00,10.00){\line(1,0){5.00}}
\put(5.00,10.00){\line(0,1){5.00}}
\put(10.00,10.00){\line(0,1){5.00}}
\put(5.00,15.00){\line(1,0){5.00}}
\linethickness{0.15mm}
\put(10.00,5.00){\line(1,0){5.00}}
\put(10.00,5.00){\line(0,1){5.00}}
\put(15.00,5.00){\line(0,1){5.00}}
\put(10.00,10.00){\line(1,0){5.00}}
\end{picture}
\end{center}}} 
& s_{M}(t)  = 3 + 2t + 1 \\ 
& 0 \r A(-4) \r A \r M \r 0 \\ 
\hline
\end{array} 
\] 

\ifx\undefined\bysame 
\newcommand{\bysame}{\leavevmode\hbox to3em{\hrulefill}\,} 
\fi

\end{document}